\documentclass[11pt,a4paper]{article}

\title{Local-global intersection homology}

\author{Jonathan Fine\relax
\thanks{203 Coldhams Lane, Cambridge, CB1 3HY, England.  
\quad E-mail: \texttt{j.fine@pmms.cam.ac.uk}}
}
\date{10 September 1997}

\newcommand\PDelta{{\bf P\!}_\Delta}
\newcommand\bfP{{\bf P}}
\newcommand\bfA{{\bf A}}
\newcommand\Abar{\bar{A}}
\newcommand\Ztilde{\widetilde{Z}}
\newcommand\Itilde{\widetilde{I}}
\newcommand\htilde{\widetilde{h}}
\newcommand\Itildehtilde{\Itilde\,\htilde}
\newcommand\Ctilde{\widetilde{C}}
\newcommand\hquest{h^?}
\newcommand\setinfty{\{\infty\}}
\newcommand\setzero{\{0\}}
\newcommand\omegaZ{\omega_Z}
\newcommand\omegaCZ{\omega_{CZ}}
\newcommand\sqdot{\bullet}
\newcommand\implies{\Rightarrow}
\newcommand\bibrule{\rule{2pc}{0.4pt}}

\begin{document}

\maketitle
\begin{abstract}\noindent
This paper defines new intersection homology groups.  The basic idea is
this.  Ordinary homology is locally trivial. Intersection homology is
not.  It may have significant local cycles.  A local-global cycle is
defined to be a family of such local cycles that is, at the same time, a
global cycle.  The motivating problem is the numerical characterisation of
the flag vectors of convex polytopes.  Central is a study of the cycles on
a cone and a cylinder, in terms of those on the base.   This leads to the
topological definition of local-global intersection homology, and a
formula for the expected Betti numbers of toric varieties.  Various
related questions are also discussed.
\end{abstract}

\section{Introduction}

This paper defines new intersection homology groups.  They record, in a
global way, local information about the singularities. They give rise to
new information, both globally and locally, and vanish on nonsingular
varieties.  Such groups are required, to obtain a satisfactory
understanding of general convex polytopes.  They also have other 
applications.

The basic idea is this.  Ordinary homology is locally trivial. 
Intersection homology is not.  It may have significant local cycles.  A
local-global cycle is a family of such local cycles that is, at the same
time, a global cycle.  The chains, that produce the homology relations
between the cycles, are to have a similar local-global nature.

The theory of toric algebraic varieties, which associates an algebraic
variety $\PDelta$ to each convex polytope $\Delta$ (provided $\Delta$ has
rational vertices) establishes a dictionary between convex polytopes and
algebraic varieties.  Convex polytopes (or, if one prefers, the associated
varieties) provide the simplest examples for these new concepts.

The basic problem is to understand general polytopes in the same way as
simple polytopes are already understood.  Suppose $\Delta$ is a simple
polytope.  Loosely speaking, this means that the associated variety
$\PDelta$ is nonsingular.  The associated homology ring $H_\bullet\Delta$
has the following properties. It is generated by the facets of $\Delta$.
It satisfies the Poincar\'e duality and strong Lefschetz theorems.  The
associated Betti numbers $h\Delta$ are a linear function of the face
vector $f\Delta$, and vice versa.  

These facts are central to Stanley's proof \cite{bib.RS.NFSP} of the
necessity of McMullen's numerical conditions \cite{bib.McM.NFSP} on the
face vectors of simple polytopes.  (An ingenious construction of Billera
and Lee \cite{bib.LB-CL.SMC} proves sufficiency.)  One would like to
understand general polytopes in a similar way.

The first results in this direction are due to Bayer and Billera
\cite{bib.MB-LB.gDS}.  They consider the flag vector, not the face vector.
For simple polytopes Poincar\'e duality represents what are known as the
Dehn-Somerville equations on the face vector.  Bayer and Billera describe
the generalised Dehn-Somerville equations on the flag vector.  They also
show that the flag vectors of $n$-dimensional polytopes span a space whose
dimension is the $(n+1)$st Fibonacci number $F_{n+1}$.

The problem of characterizing the flag vectors of general polytopes has
guided the development of local-global intersection homology.  The usual
middle perversity intersection homology theory produces $\lfloor n/2
\rfloor + 1$ independent independent linear functions of the flag vector. 
(This is the Bernstein-Khovanskii-MacPherson formula for the mpih Betti
numbers \cite{bib.JD-FL.IHNP,bib.KF.IHTV,bib.RS.GHV}.)  Clearly, more
Betti numbers are needed, to record the whole of the flag vector.  In
addition, some analogue or extension to the usual ring structure on the
homology of a nonsingular variety is required. 

The general polytope problem makes it clear that some extension of
intersection homology, and of the ring structure, is required.  Topology
by itself has failed to indicate clearly either the need for such an
extension, or its form.  (There are intersection homology groups for
non-middle perversities, and `change of perversity' groups, but these have
the same problems as ordinary homology.) Finally, there are polytopes
whose combinatorial structure is such that it cannot be realised with a
polytope that has rational vertices \cite[p94]{bib.BG.CP}.  Thus, for
general polytopes a theory that does not rely on algebraic geometry is
required.

The \emph{root problem} of this paper is as follows.  \emph{Suppose $Z$ is
a possibly singular projective algebraic variety.  In terms of the cycles
on $Z$, what are the cycles on $CZ$ and $IZ$?}  Here, $CZ$ is the
projective cone on $Z$, while $IZ$ is the product of $Z$ with $\bfP_1$. 
The answer depends on what one understands a cycle to be, or in other
words on some perhaps implicit choice of a homology theory. Suppose this
question has been answered.  One will then have a wide range of examples.
These will determine the corresponding definition of a cycle, in the same
way that a number of points will determine a plane.

These examples will also determine a linear function $h\Delta$ of the flag
vector $f\Delta$ of convex polytopes.  This is because of the following. 
There are operators $I$ and $C$ on polytopes, analogous to the $I$ and $C$
operators on varieties.  An \emph{$IC$ polytope} is any polytope that can
be obtained by repeatedly applying $I$ and $C$ to the point polytope.  Any
polytope flag vector can be written as a linear combination of the flag
vectors of $IC$ polytopes.  The examples determine $h\Delta$ on the $IC$
polytopes.

This paper is organised as follows.  First (\S2) notation and definitions
are established, and some basic results stated.  Next (\S3) the root
problem is discussed and a solution presented.  This serves to motivate
the definition of the (extended) $h$-vector $h\Delta$ of convex polytopes
(\S4), and the topological definition of local-global intersection homology
(\S5).  To finish (\S6), there is a summary, and a discussion of related
questions.

This paper considers the topological and combinatorial aspects of
local-global homology.  There are others, to be presented elsewhere.  The
\emph{linear algebra} \cite{bib.JF.CPLA} allows $h\Delta$ to be
interpreted as the outcome of a `vector weighted inclusion-exclusion'
construction.  The \emph{intersection theory} \cite{bib.JF.IHRS} provides
a structure that reduces, in the simple or nonsingular case, to the ring
structure on ordinary homology.

To date, the theory of local-global intersection homology consists of a
series of definitions appropriate for each of the four aspects, together
with examples and special cases, and various linking results.  Much remains
to be done, to fill in the `convex hull' of the four aspects.

This paper has been written to be accessible to those who are unfamiliar
with perhaps one or both of intersection homology and the combinatorics of
convex polytopes.  The reader who is in a hurry can find a summary in the
final section.  Formulae (\ref{eqn.Itilde})--(\ref{eqn.Abar-A}) define the
extended $h$-vector.  The topological definition is in \S5.  Text in
parentheses (except for short comments) can be omitted on a first
reading.  The reader who is having difficulties should first understand
the mpih part of the theory (i.e.~ignore terms involving any of $A$,
$\Abar$ or $\{k\}$).

\section{Preliminaries}

This section introduces notation and conventions.  It also states results
to be used later.  This material is organised into six topics, namely
language and conventions, cones and cylinders, local homology, the strong
Lefschetz theorem, polytope flag vectors, and the index set for $h\Delta$.

First, language and conventions.  Ordinary homology fails to have suitable
properties, and so the word `homology' when used without qualification
will refer either to middle perversity intersection homology (mpih), or
some local-global variant thereof.  The abbreviation mpih will always
refer to the usual intersection homology (as in \cite{bib.MG-RDM.IH}),
with of course middle perversity.  Each local-global homology group has an
order, usually denoted by $r$. The usual mpih groups are order zero
local-global groups.  The higher order groups will be called
\emph{strictly local-global}.  Unless otherwise stated, homology will
always be with rational (or real) coefficients.

The analogy between convex polytopes and algebraic varieties is very
important, particularly in \S3 and \S4.  Throughout $\Delta$ will be a
convex polytope of dimension~$n$, and $Z$ a (projective) algebraic
variety, also of dimension~$n$.  When $\Delta$ has rational vertices a
projective toric variety $\PDelta$ (of dimension~$n$) can be constructed
(as in, say, \cite{bib.VD.GTV}). If $\Delta$ is the $n$-simplex $\sigma_n$
then $\PDelta$ is projective $n$-space $\bfP_n$.

To strengthen the analogy, for toric varieties homology will be indexed by
the complex dimension (half the normally used real dimension).  The mpih
Betti numbers of $\PDelta$ are zero in the odd (real) dimensions, and the
same is expected to hold for the local-global extension.  Thus, this
indexing convention amounts to ignoring the homology groups that are
expected in any case to be zero.

The concept of a cone is one of the most important in this paper.  In
fact, the same word will be used for three closely related constructions,
that apply respectively to topological spaces, projective algebraic
varieties, and convex polytopes.

Suppose $B$ is a topological space.  The \emph{cylinder} $IB$ on $B$ is
the product $[0,1]\times B$ of $B$ with the interval $I=[0,1]$, equipped
with the product topology.  If $p=(\mu,l)$ is a point of $IB$ and
$\lambda\in I$ is a scalar then $\lambda p = (\lambda\mu,l)$ is also a
point on $IB$.  The \emph{cone} $CB$ is the cylinder $IB$, with
$\setzero\times B$ identified (collapsed) to a single point, the
\emph{apex} of the cone.  In \S5, local-global cycles will be described as
global cycles that can be collapsed in some specified way.  The locus
$\{1\}\times B$ is called the \emph{base} of $CB$.  There is a
$\lambda$-action on $CB$ also.

Now suppose $Z\subset \bfP_N$ is a projective algebraic variety.  The
(projective) \emph{cone} $CZ\subset \bfP_{N+1}$ with \emph{base} $Z$ is
constructed as follows.  Each point $p\in\bfP_N$ represents a line $l_p$
through the origin in $\bfA^{N+1}$.  A point $v$ lies in the \emph{affine
cone} $\Ztilde \subset \bfA^{N+1}$ just in case it lies on some $l_z$,
with $z\in Z$.  A `hyperplane at infinity' can be added to $\bfA^{N+1}$,
to produce $\bfP_{N+1}$.  This hyperplane is a copy of $\bfP_N$.  The cone
$CZ$ is the closure of $\Ztilde$.  The base of $CZ$ is the copy
$\setinfty\times Z$ of $Z$ that lies on the $\bfP_N$ at infinity.  The
origin of $\bfA^{N+1}$ is the \emph{apex} of the cone.  

The interaction between the cone structure and relations among cycles is
central to \S3.  The complex numbers act by multiplication on the `finite
part' $\Ztilde$ of $CZ$.  Thus, if $\eta$ is a cycle on $CZ$, lying
entirely on $\Ztilde$, it can as be `coned away' by the
\emph{$\lambda$-cone} $C_\lambda\eta$, where $\lambda$ ranges over
$[0,1]$.  This is a chain whose boundary is $\eta$ (unless $\eta$ has
dimension zero).  Now suppose $\eta$ on $CZ$ avoids the apex.  In this case
each point of $\eta$ lies on a unique line through the apex, and so there
is a boundary that `moves' $\eta$ to an equivalent cycle $\eta_\infty$
that lies entirely on the base $Z$ of $CZ$.  Finally, suppose that $\eta$
is a cycle lying on the base $Z$ of $CZ$.  Each point of $\eta$ determines
a line in $CZ$, and so $\eta$ determines a cycle $C\eta$ on $CZ$. However,
it may not be possible to find an $\eta'$ lying entirely on $\Ztilde$,
that is equivalent to $\eta$.  

(The reason for this is subtle.  If it were always possible, then it would
be possible to `cone away' the cycle due to the hyperplane at infinity in
$\bfP_{n+1}$ (the cone on $\bfP_n$).  But this cycle is not homologous to
zero.  Although $\eta$ can locally be moved away from the base of $CZ$, in
a manner that is unique up to `phase', it may not be possible to get all
the phases to match up.)

The \emph{cylinder} $IZ$ is the product of $Z$ with $\bfP_1$, which via
the Segre embedding is to be thought of as a subvariety of $\bfP_{2N+1}$. 
The variety $Z$ is the \emph{base} of the cylinder.  If $Z$ is
nonsingular, then so is $IZ$, whereas $CZ$ will in general have a
singularity at its apex.

Analogous operators $I$ and $C$ can be defined for convex polytopes.  If
$\Delta$ is a convex polytope then the \emph{cone} (or \emph{pyramid})
with \emph{base} $\Delta$ is the convex hull of $\Delta$ with a point (the
\emph{apex}) that does not lie in the affine span of $\Delta$.  Similarly,
the \emph{cylinder} (or \emph{prism}) $I\Delta$ with \emph{base} $\Delta$
is the Cartesian product $[0,1]\times\Delta$ of $\Delta$ with an interval
$I=[0,1]$.  These operators respect the dictionary between convex
polytopes and toric algebraic varieties.

The symbol `$\sqdot$' will be used to denote both the projective variety
$\bfP_0$ (a single point), and the single point convex polytope.  Thus,
$ICC\sqdot$ can denote either $\bfP_1\times\bfP_2$ or a triangular prism.
An \emph{$IC$ polytope} is one obtained by successively applying $I$ and
$C$ to the point polytope, and similarly for an \emph{$IC$ variety}.  For
every word in $I$ and $C$, the latter is the toric variety associated to
the former.  Sometimes the two concepts will be identified.

In both cases one can also define the \emph{join} of two objects. Suppose
that $Z_1$ and $Z_2$ are subvarieties of some $\bfP_N$, and that their
affine linear spans are disjoint.  In that case, their \emph{join}
consists of all points that lie on some line $l(z_1,z_2)$ that joins a
point $z_1\in Z_1$ to another $z_2\in Z_2$.  Similarly, if $\Delta_1$ and
$\Delta_2$ have disjoint affine linear spans, then their \emph{join} is
the convex hull of $\Delta_1\cup\Delta_2$.  For both polytopes and
varieties, a cone is the join of the base to the apex.  One could also
join an object not to a point but to a projective line (respectively, an
interval).  This is the same as forming the cone on the cone.  It will
have an \emph{apex line} (resp.~\emph{apex edge}) rather than an apex.

Intersection homology differs from ordinary homology in that for it
nontrivial local cycles can exist.  If $s$ is a point on a complex
algebraic variety $Z$ any sufficiently small ball centered at $s$ is
homeomorphic to the (topological) cone $CL_s$ on something.  That
something, which does not depend on the sufficiently small ball, is the
\emph{link} $L_s$ at $s$. Now suppose $\eta$ is a cycle on $CL_s$.  The
$\lambda$-action on a cone can then be applied to $\eta$, to produce a
chain $C_\lambda\eta$, whose boundary is $\eta$.  (Strictly speaking, this
is true only if $\dim \eta > 0$.)  Local cycles of dimension zero are
trivial, and will not be counted by $h\Delta$.  Ordinary homology allows
this \emph{coning away} of local cycles. The perversity conditions of
intersection homology however can be used to prohibit the use of
$C_\lambda\eta$ to generate a boundary.

The local (intersection) homology groups can be defined as follows.  A
\emph{local cycle} $\eta$ at $s$ consists of a cycle $\eta_U$ for any
sufficiently small open set $U$ containing $s$, such that if $U'\subset
U$, one has on $U$ that $\eta_{U'}$ and $\eta_U$ are homologous. 
Similarly, a \emph{local boundary} $\xi$ at $s$ consists of a chain
$\xi_U$ on each sufficiently small open set $U$, whose boundary $\eta$
(the system $\eta_U=d\xi_U$) is a local cycle.  This definition avoids use
of the cone structure.  In \S4, another definition will be given.

If $Z$ is a complex algebraic variety then it can be decomposed into a
disjoint union of \emph{strata} $S_i$, where each $S_i$ is either empty or
has complex dimension~$i$, and along $S_i$ the local topology of $Z$ is
locally constant.  From this it follows that the local homology groups are
also locally constant along $S_i$, and so form what is known as a
\emph{local system}.  This concept is used only in \S5.  However, it is
closely related to an example of local-global homology.

(This paragraph and the next can be omitted on a first reading.)  Suppose
that $i<j$ and that the stratum $S_i$ is in the closure of $S_j$.  More
particularly, suppose that $\gamma:[0,1]\to Z$ is a path, with
$\gamma(0)\in S_i$ and $\gamma(\lambda)\in S_j$ otherwise.  Now let
$\eta_1$ be a local cycle at $\gamma(1)$.  By local constancy, it can be
moved along $\gamma$ until it is very close to $\gamma(0)$.  At this point
the translate $\eta_\lambda$ of $\eta_1$ can be thought of as a local
cycle at $\gamma(0)$ on $S_i$.  In other words, each path from $S_j$ to
$S_i$ (with $j>i$) transfers local cycles from $S_j$ to $S_i$.  Note that
the reverse process will not in general be possible.  For example, if
$S_0$ is an isolated singularity, then a local cycle $\eta$ at $S_0$
cannot be moved away from $S_0$.

Now consider $H_0(S_i,L_i)$, where $L_i$ is the local system formed from
the local homology groups along $S_i$.  A cycle $\eta\in H_0(S_i,L_i)$ is
a formal sum of local homology cycles (about points on $S_i$) subject to
the equivalence due to motion along paths.  As already described, these
groups can be `glued together' (certain elements identified) for
different values of $i$.  Provided one uses all $S_i$ whose dimension is
at least some value value $j$, the result is independent of the
stratification.  (This is left to the reader.  Then main point is that new
strata have real codimension at least two, and so existing paths can be
altered to avoid new strata.)  These groups are examples of local-global
homology.

The \emph{strong Lefschetz theorem} is one of the central results in the
homology of nonsingular algebraic varieties.  It was stated by Lefschetz
in 1924, but his proof was not satisfactory.  The first proof is due to
Hodge (1933--6, see \cite[p117]{bib.WVDH.TAHI}).  It also follows from
Deligne's proof of the Weil conjectures \cite{bib.NK.DPRH}.  Strong
Lefschetz also holds for middle perversity intersection homology.  Here,
Deligne's proof is the only method known.  For more background see
\cite{bib.SK.DIHT}.  The infinitesimal form of Minkowski's facet area
theorem for polytopes~\cite[p.332]{bib.HM.ALKP,bib.BG.CP} is a special
case of both strong Lefschetz and the Riemann-Hodge inequalities.  
This seems not to have been noticed before.

Suppose $Z\subset\bfP_N$ is a projective algebraic variety.  For
convenience, complex dimension will be used to index its homology groups
$H_iZ$.  If $i+j=n$ (the dimension of $Z$), then by Poincar\'e duality
$H_iZ$ and $H_jZ$ have the same dimension, for they are dual vector
spaces.  The embedding $Z\subset \bfP_N$ determines a \emph{hyperplane
class} $\omega=\omegaZ$ in $H_{n-1}Z$ with the following properties. 
First, the cap product $\omega\frown\eta$ is defined for any homology
class $\eta$ on $Z$.  This operation lowers degree by one.  Now assume
$i<j$ and also $i+j=n$.  The \emph{strong Lefschetz theorem} asserts that
the map
\[
    \omega^{j-i}: H_jZ \to H_iZ
\]
is an isomorphism.

This result provides a decomposition of $H_0Z$.  Suppose that the above
isomorphism takes $\eta$ to $\eta'=\omega^{i-j}\eta$.  Say that $\eta$ is
\emph{primitive} and $\eta'$ is \emph{coprimitive}, if $\omega\eta'$ is
zero.  It is a standard result, that $H_\bullet Z$ is the direct sum of
$\omega^iP_jZ$, where $P_jZ\subseteq H_jZ$ are the primitive classes, and
$i+j\leq n$.  The Lefschetz isomorphism allows an `inverse' $\omega^{-1}$
to $\omega$ to be defined.  Define $\omega^{-1}$ to be the result of first
applying the inverse of the Lefschetz isomorphism, then $\omega$, and then
the Lefschetz isomorphism.  It has degree $-1$.  (The Riemann-Hodge
inequality is that on $P_jZ$ the quadratic form
$\eta\frown\omega^{j-i}\frown\eta$ is negative definite.)

The primitive and coprimitive cycles have a special r\^ole in the study of
the root problem, namely the cycles on $CZ$ in terms of those on $Z$. They
can occur only in certain dimensions, for which it is useful to have
special adjectives.  Say that a cycle $\eta$ on $Z$ is \emph{upper}
(respectively \emph{strictly upper}) if its dimension is at least
(resp.~more than) half that of $Z$.  Similarly, at most (resp.~less than)
define \emph{lower} (resp.~\emph{stricly lower}).  Primitives and
coprimitives occur in the upper and lower dimensions respectively.  A
cycle dimension that is not lower is strictly upper, and vice versa.  The
\emph{middle dimension} is both upper and lower.

The hyperplane class $\omegaZ$ on $Z$ can be represented as a Weil divisor
(formal sum of codimension one subvarieties) on $Z$, namely the hyperplane
section.  For use in \S3, note that the hyperplane class $\omegaCZ$ on a
cone can be represented either as the cone $C\omegaZ$ on the class of the
base, or as the base $Z$ of the cone (by intersecting $CZ \subset
\bfP_{N+1}$ with the $\bfP_N$ at infinity).

If $\Delta$ is a \emph{simple} convex polytope (this means that at each
vertex there are $n=\dim\Delta$ edges) then $\PDelta$ behaves like a
nonsingular algebraic variety, so far as its homology (with rational
coefficients) is concerned.  Its Betti numbers $h_i\Delta=h_i\PDelta$ are
then a linear function of the \emph{face vector}
$f=f\Delta=(f_0,f_1,\ldots,f_n)$, where $f_i$ is the number of
$i$-dimensional faces on $\Delta$.  In fact, if one writes $f(x,y)=\sum
f_ix^iy^{n-i}$, and $h(x,y)$ similarly, then the equation
$h(x,x+y)=f(x,y)$ expresses the relation between $f$ and $h$.

If $\Delta$ is a general convex polytope, then flags should be counted.  A
\emph{flag} is a sequence
\[
    \delta = ( \delta_1 \subset \delta_2 \subset
        \ldots \subset \delta_r \subset \Delta )
\]
of faces, each stricly contained in the next.  Its \emph{dimension vector}
(or \emph{dimension} for short) is the sequence
\[
    d = ( d_1 < d_2 < \ldots < d_r < n )
\]
of the dimensions $d_i$ of its terms $\delta_i$.  Altogether, there are
$2^n$ possible flag dimensions.  The component $f_d\Delta$ of the
\emph{flag vector} $f=f\Delta$ of $\Delta$ counts how many flags there on
$\Delta$, whose dimension is $d$.  (If $\Delta$ is simple, the flag vector
is a linear function of the face vector, and so contains no new
information.)

For simple polytopes the \emph{Dehn-Somerville} equations state that
$h(x,y)$ is equal to $h(y,x)$, or that the $h$-vector is
\emph{palindromic}.  (It is analogous to Poincar\'e duality.)  For general
polytopes the \emph{generalised Dehn-Somerville (gDS) equations}
\cite{bib.MB-LB.gDS} imply that $f\Delta$ has the Fibonacci number
$F_{n+1}$ linearly independent components.  A similarly elegant
interpretation of these equations is lacking.

For flag vectors the \emph{$IC$ equation} \cite{bib.JF.MVIC}
\[
    (I-C)C \> I \> = \> I \> (I-C)C
\]
holds, in the following sense.  Apply both sides to a convex polytope
$\Delta$, to obtain convex polytopes $ICI\Delta$ etc.  The corresponding
equation then holds among the flag vectors of these polytopes. The flag
vectors of the $IC$ polytopes span all polytope flag vectors, and those
than contain neither `$II$' nor `$I\sqdot$' form a basis.

It follows that if linear operators $\Itilde$ and $\Ctilde$ are given that
satisfy the $IC$ equation, together with an initial value
$\htilde(\sqdot)$ for which $\Itildehtilde(\sqdot)=\Ctilde\htilde(\sqdot)$,
then there is a unique linear function $\htilde$ on polytope flag vectors,
for which the equations
\[
    \htilde(I\Delta) = \Itilde h(\Delta) \>;\quad
    \htilde(C\Delta) = \Ctilde h(\Delta) \>;
\]
are satisfied.  This is used in \S4, to define the extended $h$-vector.

Finally, note that the polytope flag vectors span a proper subspace of the
span of all flag vectors.  To provide a linear function on this subspace
is not the same as to provide such on the larger space.  Conversely,
different linear functions on the larger space can agree on the subspace.
Related to this is the idea that equivalent homology theories (on
projective varieties) can be given different definitions, and that
different triangulations can be found for a given space.

The last topic is the \emph{index set}.  The extended $h$-vector
$h\Delta$, to be defined in \S4, will be a formal sum of terms, of a
particular type.  Although the terms to be used will arise in a natural
way, it is convenient to gather in one place a description of them.  In
some sense, their structure is a result of a geometric requirement (that
the Betti numbers be organised into sequences, whose length grows with the
dimension $n$ of the polytope) and a combinatorial requirement (the
Fibonacci numbers).

First, expressions such as $(a,b,c)$ will stand for the homogeneous
polynomial $ax^2 +bxy + cy^2$ in commuting variables $x$ and $y$.  To save
space, commas will where possible be omitted.  Thus, $(10)=x$, $(01)=y$,
$(11)=x+y$ and $(1)=1$.  Similarly, the expression $[abc]$ (short for
$[a,b,c]$) will stand for $aX^2 +bXY +cY^2$, where $X$ and $Y$ are a
different pair of commuting variables.  Clearly, $[1]=(1)=1$.  Each of
$x$, $y$, $X$ and $Y$ will have degree one.

Roughly speaking, to each coprimitive cycle on the base of a cone (or in
the link along a face) the symbol $\{k\}$ will correspond, where $k$ is
the dimension of the cycle.  This symbol also corresponds to a local
$k$-cycle.  Because in the present context such can occur only in the
strictly lower dimensions, $\{k\}$ will be given degree $2k+1$.  The
symbol $\{0\}$ will not be used.  This corresponds to treating the class of
a point as a trivial local cycle. The difference $b-a$, which `counts'
dimension~$1$ coprimitives, will be denoted by $b'$.  Similarly, $c'=c-b$,
and so on up to halfway.  In addition, a `padding' symbol $A$ (or $\Abar$)
is required. It has degree one.

The extended $h$-vector $h\Delta$ will be a sum of terms of the form
$x^iy^jW$, where $W$ is a word in $A$ and $\{k\}$.  Each such term will
have degree equal to the dimension of $\Delta$.  The word $W$ is allowed
to be empty.  This corresponds to the mpih part of $h\Delta$.  The last
symbol in $W$ is not to be an $A$.  (This can be achieved by supposing
that there is a terminating symbol `$\sqdot$' at the end of each word, and
setting $A\sqdot$ equal to zero.) An auxiliary vector $\htilde\Delta$ is
used in \S4.  It is a sum of $X^iY^jW$ terms, where $W$ is a word in
$\Abar$ and $\{k\}$.  Its terms are otherwise the same as those of
$h\Delta$.

The numerology of $x^iy^jW$ is interesting.  Recall that $h\Delta$ has
$F_{n+1}$ independent components.  There are $F_{n+2}$ terms satisfying
the above conditions, whose degree is $n$.  Of these $F_{n+1}$ satisfy
$i\leq j$, (and they correspond to a maximal set of independent components
in $h\Delta$).  Thus, $F_n$ ($=F_{n+2}-F_{n+1}$) terms satisfy $i>j$. 
Similarly, $F_n$ terms satisfy $j>i$.  Thus, $F_{n-1}$ terms satisfy
$i=j$.  The number of words $W$ of degree at most $n$ is $F_n$, for $n\geq
1$.  (These results are not used, and so are stated without proof.)

Also, the equation
\[
    1 + F_1 + F_2 + \ldots + F_n = F_{n+2}
\]
can be interpeted as follows.  Define $f_{(i)}\Delta$ to be the sum of the
flag vectors of the $i$-faces (or $i$-links if one prefers) of $\Delta$. 
It follows from the flag vector concept that $f\Delta$ and $f_{(\bullet)}=
(1, f_{(1)}, f_{(2)}, \ldots, f_{(n-1)})$ are linear functions of each
other.  (The `$1$' corresponds to the `empty' face, or to $\Delta$
itself.)  Each $f_{(i)}$ has $F_{i+1}$ independent components, and so
$f_{(\bullet)}\Delta$ (and hence $f\Delta$) has $F_{n+2}$ components whose
dependence does not follow from the gDS equations on the faces.  This
helps justify $F_{n+2}$ as the number of components in $h\Delta$, for one
wishes $h\Delta$ to permit an elegant expression of the generalised
Dehn-Somerville equations.

\section{Cycles on cones and cylinders}

This section describes the mpih and local-global cycles on a cone and a
cylinder in terms of those on the base.  First, the mpih cycles are
constructed and described.  The local-global cycles are then a variant of
the mpih cycles.  They make use of information, that mpih ignores.

Suppose that $\eta$ is a cycle on $Z$.  Later, this statement will acquire
a richer meaning, but for now suppose that $\eta$ is a formal sum of
embedded simplices, whose boundary is zero.  The cycle $\eta$ on $Z$
determines three cycles on $IZ$, which can be denoted by
$\setzero\times\eta$, $\setinfty\times\eta$, and $I\eta$.  The first two,
which are equivalent, arise from the two `poles' $0$ and $\infty$ on
$\bfP_1$, each of which determines an embedding of $Z$ in $IZ$.  The
third, $I\eta$, is the product of $\eta$ with $\bfP_1$.  Similarly, on
$CZ$ one will have $\setinfty\times\eta$ and $C\eta$.  (There is also the
apex, which will not be needed.)

Relations, as well as cycles, must be considered.  On $\bfP_1$ let
$I_\lambda$ denote the chain that is a path from $0$ to $\infty$. 
Similarly, let $I_\lambda\eta$ denote the chain on $IZ$, whose boundary is
$\setinfty\times\eta -\setzero\times\eta$.

The relations on a cone are more complicated.  As noted in \S2, in general
it is not possible to cone a cycle $\eta$ on the base to form a relation
$C_\lambda\eta$, whose boundary is $\eta$. However, if a cycle $\eta$ on
$CZ$ is equivalent to an $\eta'$ that does not meet the base $Z$, then
$\eta$ can (via $\eta'$) be coned away to produce a relation
$C_\lambda\eta$.  Note that if such an $\eta'$ can be found, then the cap
product of $\eta$ with the base (if defined) will be equivalent to zero,
for $\eta'$ does not meet the base.

When no restrictions are places on the cycles and relations, ordinary
homology is the result.  Based on the preceeding discussion, one might
expect the ordinary homology of $IZ$ to be the tensor product of that of
$Z$ with that of $\bfP_1$ (the K\"unneth formula), while for $CZ$ one
might expect the `cone' on the ordinary homology of the base.  By this is
meant the base homology raised by one in degree, with the class of a point
appended in degree zero.  There are similar expected formulae for the
ordinary homology Betti numbers.  However, ordinary homology does not in
general satisfy Poincar\'e duality and strong Lefschetz.  Also, its Betti
numbers are not a linear function of the flag vector \cite{bib.McC.HTV}. 
This is discussed further in \S5.

Now consider middle perversity intersection homology, or more precisely, a
theory that satisfies the Poincar\'e and Lefschetz theorems.  These
properties, particularly strong Lefschetz, will leave one with little
choice as to what the cycles on $CZ$ and $IZ$ are, and hence lead to the
usual middle perversity conditions on the cycles.

The task is to control the cycles and relations, so that the $I$ and $C$
operators preserve the Poincar\'e and Lefschetz properties.  For $I$ the
usual K\"unneth formula will do this, a result that is left to the reader.
For $C$, more care is needed.  


Suppose $\setinfty\times\eta$ is a cycle on $CZ$.  Now use the $C\omegaZ$
form of the hyperplane class.  Clearly, one will have
$C\omegaZ\frown(\setinfty\times\eta) =
\setinfty\times(\omegaZ\frown\eta)$ as the hyperplane action.  Similarly,
if $C\eta$ is a cycle on $CZ$, use the `base' form $Z$ of the hyperplane
class to obtain
\[
    \omegaCZ \frown C\eta \sim Z \frown C\eta = \setinfty\times\eta
\]
as the hyperplane action.  Now suppose that $\eta$ is a primitive
$j$-cycle on $C$, by virtue of $\omegaZ^i\eta \sim 0$.  The relation
\[
    \omegaCZ^{i+1} C\eta \sim 0
\]
follows from the above.  Thus, \emph{$C\eta$ is a primitive on $CS$,
whenever $\eta$ is a primitive on $Z$}.  Moreover, $\setinfty\times\eta$
is equal to $\omegaCZ\frown C\eta$, and so cannot be primitive.  All this
assumes that $C\eta$ is allowed as a cycle, when $\eta$ on $Z$ primitive.

These properties (K\"unneth and the coning of primitives) suffice to
determine the homology of $IZ$ and $CZ$ respectively, in terms of that of
$Z$.  The task now is to express these groups in terms of topological
cycles and relations.  First consider $CZ$.  By assumption, if $\eta$ is
primitive on $Z$, then $C\eta$ is permitted on $CZ$ (and is there
primitive).  Primitive is not a topological concept; it depends on the
projective embedding $Z\subset\bfP_N$.  However, the primitive cycles have
upper dimension, and that is a topological notion.  Thus, permit
$\xi=C\eta$ as a cycle on $CZ$ whenever $\eta$ on $Z$ is upper, or in
other words when $\xi$ on $CZ$ is strictly upper.

Now suppose that $\eta$ on $Z$ is coprimitive.  It follows at once that
\[
    \setinfty\times\eta \frown \setinfty\times Z \sim 0 \>,
\]
and so there is nothing in the homology of $CZ$, that prevents
$\setinfty\times\eta$ being moved away from the base $\setinfty\times Z$. 
Suppose that this can be done, to produce $\eta' \sim
\setinfty\times\eta$.  As noted in \S2, the $\lambda$-coning operation can
be applied to $\eta$, via $\eta'$, to produce $C_\lambda\eta$.  If $\eta$
on $Z$ is nonzero, then on $CZ$ it is also by assumption nonzero.  Thus,
the `coning away to a point' $C_\lambda\eta$ cannot be permitted (except
perhaps if $\dim \eta =0$).  As before, even though coprimitive is not a
topological notion, the lower range of dimensions is.  This leads to
$\xi=C_\lambda\eta$ on $CZ$ being prohibited as a chain, whenever $\xi$ is
lower.

These two examples (primitive and coprimitive) establish the middle
dimension as the cut-off point for cycles and chains being permitted or
prohibited respectively.  This applies to how they meet the $0$-strata. 
To obtain the remainder of the middle perversity conditions, study the
cycles on $IZ$ due to $\eta$ on $Z$,
where now $\eta$ is a cycle that satisfies the conditions that
are already known.  In this way, the rest can be built up, to produce the
already known middle perversity intersection homology conditions on cycles
and chains.  This is left to the reader.

The assumption, that if $\eta$ on $Z$ is coprimitive, then
$\setinfty\times\eta$ is equivalent to an $\eta'$ that does not meet the
base $Z$ of the cone, is quite strong.  Previously, it was assumed that
this might happen from time to time, and so the associated coning aways
were prohibited.  That there are such prohibited coning aways is a
topological property, which is local to the apex of $CZ$.  Local-global
homology will count such `coning aways', but calls them local-global cycles. 
It will as a heuristic principle be assumed that any coprimitive can be
moved to avoid the base, and so be coned away.   Such assumptions support
the calculation in \S4 of the expected values $h\Delta$ of the
local-global Betti numbers.

It is time to take stock.  Recall that the purpose of this section is to
describe the cycles on $IZ$ and $CZ$ in terms of those on $Z$.  When no
restrictions are imposed, ordinary homology is the result.  Requiring
Poincar\'e duality and strong Lefschetz produces middle perversity
intersection homology.  For mpih the cycles on $CZ$ are all of the form
$C\eta$ or $\setinfty\times\eta$, for $\eta$ a cycle on $Z$.  For the
local-global extension, one also has the `coning away' or
\emph{local-global cycle} $C_\lambda\eta$, for $\eta$ any coprimitive
cycle on the base $Z$.  (The apex of the cone $Z$ is also called the
\emph{apex} of the local-global cycle $C_\lambda\eta$.)  Now assume that
on $Z$ itself there is such a local-global cycle.  Further local-global
cycles may arise on $CZ$ and $IZ$, as a result of this local-global cycle
on $Z$.

As in the mpih case, on $IZ$ it will be assumed that the K\"unneth formula
continues to hold.  In other words, any cycle on $IZ$ can be expressed
using $\setzero\times\eta$ (or $\setinfty\times\eta$) and $I\eta$, where
$\eta$ ranges over the cycles on $Z$.  If $\eta$ is a local-global cycle,
say $C_\lambda\xi$, then $I\eta$ is a new kind of local-global cycle. (Its
apex locus is $I$ applied to that of $\eta$.)  If $\eta$ is thought of as
a local cycle, then $I\eta$ is a family of local cycles.  In addition,
note that by K\"unneth $\setzero\times\eta$ and $\setinfty\times\eta$ are
equivalent cycles on $IZ$, and so local-global cycles must on occasion be
allowed to move along the singular locus.  (In fact, the rule will be that
if they can move, then they are allowed to move.)

At this point it is possible to give some examples.  In dimensions $0$,
$1$ and $2$ one has
\[
    h(\sqdot) = (1) \>\>; \quad
    h(C\sqdot) = h(I\sqdot) = (11) \>\>; \quad
    h(CC\sqdot) = (111) \>\>; \quad
    h(IC\sqdot) = (121) \>\>; 
\]
of course.  In dimension $3$ one has
\[
    h(CCC\sqdot) = (1111) \>\>;\quad
    h(ICC\sqdot) = (1221) \>\>;\quad
    h(IIC\sqdot) = (1331) \>\>;
\]
as the nonsingular (or simple) examples, while
\[
    h(CIC\sqdot) = (1221) + (1)\{1\} \>\>;
\]
is the only singular example.  Here, $\{1\}$ counts the local-global
cycles on $CZ$ due to the only nontrivial coprimitive on $Z$, where $Z$ is
$\bfP_1\times\bfP_1$ (or a square).

In dimension $4$ something new happens.  The simple cases
\[\begin{array}{ll}
    h(CCCC\sqdot) = (11111) \>\>;\quad &
    h(ICCC\sqdot) = (12221) \>\>;\\
    h(IICC\sqdot) = (13431) \>\>;\quad &
    h(IIIC\sqdot) = (14641) \>\>;\\
\end{array}\]
are just as before.  The cones on the simple dimension $3$ examples come
next.  They are
\[
    h(CICC\sqdot) = (12221) + (1)\{1\}A \>\>;\quad
    h(CIIC\sqdot) = (13331) + (2)\{1\}A \>\>;
\]
where as before $\{1\}A$ counts the nontrivial coprimitives on the base. 
The remaining examples are $I$ and $C$ applied to $CIC\sqdot$, the only
non-simple dimension $3$ example.

The polytopes (or varieties) $ICIC\sqdot$ and $CCIC\sqdot$ are more similar
than they might at first sight appear.  The polytope $CIC\sqdot$ has an
apex, and so $ICIC\sqdot$ has an \emph{apex edge}.  Now consider
$CCIC\sqdot$.  From one point of view, this has two apexes, namely the
apex of its base $CIC\sqdot$, and the apex of $CCIC\sqdot$ itself. 
However, $CCIC\sqdot$ is also the join of $IC\sqdot$ with an interval, and
so there is no geometric way of distinguishing its two apexes.  In other
words, like $ICIC\sqdot$, it too has an apex edge.  Combinatorially, the
two polytopes are the same along their respective apex edges.  (This fact
is at the heart of the $IC$ equation.)

The mpih parts of $h(ICIC\sqdot)$ and $h(CCIC\sqdot)$ are $(13431)$ and
$(12221)$ respectively.  The remaining contribution comes from the
strictly local-global cycles along the apex edge.  Clearly, if $\eta$ is
the cycle that contributes $(1)\{1\}$ to $CIC\sqdot$, then $I\eta$ and
$\setzero\times\eta$ will contribute $(01)\{1\}$ and $(10)\{1\}$
respectively to $ICIC\sqdot$. However, on $CIC\sqdot$ there is a single
coprimitive cycle (it has dimension one), and so on $CCIC\sqdot$ there
will be a local-global of type $\{1\}A$.  Now note that if $h\Delta$ is to
be a linear function function of the flag vector, the non-mpih parts of
the $h$-vectors of the two polytopes should be the same.  To achieve this,
the values
\[\begin{array}{l}
    h(ICIC\sqdot) = (13431) + (11)\{1\} + \{1\}A \>\>\;\\
    h(CCIC\sqdot) = (12221) + (11)\{1\} + \{1\}A \>\>\;
\end{array}\]
will be postulated. 

This is to take any strictly local-global contribution that can occur for
either $ICIC\sqdot$ or $CICC\sqdot$, and to insist that it can occur in
the other.  This forces $\setzero\times\eta$ on $ICIC\sqdot$ to contribute
not only $(10)\{1\}$ as already noted, but also $(1)\{1\}A$. The cycle
$\setinfty\times\eta$ on $CZ$ will make a similar contribution to
$hCCIC\sqdot$.  Also, some sort of coning $C\eta$ of the cycle $\eta$ must
be allowed, to obtain on $CCIC\sqdot$ a contribution of $(10)\{1\}$. This
discussion is an example of how topology and combinatorics work together
to determine the structure of the theory of local-global homology.

Now suppose that $\eta$ is a cycle (possibly of local-global type) on $Z$.
Already, the associated cycles on $IZ$ have been described.  The task now
is to determine and describe the associated cycles on $CZ$.  There are
three basic possibilities.  First, one can form $\setinfty\times\eta$,
which is a cycle lying on the base $Z$ of $CZ$.  Second, one can cone
$\eta$ to form $C\eta$.  Sometimes, as in the lower dimensions of mpih,
this cycle is not needed.  Finally, if $\setinfty\times\eta$ can be moved
to an equivalent cycle $\eta'$, that does not meet the base $Z$ of $CZ$,
one can form the `coning away' $C_\lambda\eta$.  These possibilities will
be considered, one at a time.

First, the cycle $\setinfty\times\eta$ will always be admitted.  There are
no conditions imposed on $\eta$.  One reason for this is that about their
respective bases, the cylinder and the cone are combinatorially the same,
and so that which is permitted for the one should be permitted for the
other.  But for the cylinder, the K\"unneth principle causes
$\setinfty\times\eta$ to be admitted.  As in the cylinder, this cycle may
contribute to several distinct local-global homology groups.

Next consider $C\eta$.  The example of $CCIC\sqdot$ shows that this case
requires more thought.  As already noted, the cycle $\eta$ on $CIC\sqdot$
contributes $(11)\{1\}$ plus $(1)\{1\}A$ to $CCIC\sqdot$.  Clearly,
$\setinfty\times\eta$ contributes $(10)\{1\}$ and $(1)\{1\}A$.  The
remainder, $(01)\{1\}$ will have to come from $C\eta$.  Thus, at least in
this case, $C\eta$ must be allowed.  But this seems to contradict the mpih
case, where no use of $C\eta$ was made in the lower dimensions, and where
any $C_\lambda\eta$ `coning away' was explicitly prohibited.  However, it
is possible to harmonise the two cases.  This involves looking again at
the mpih situation.

Suppose that $\eta$ is an mpih cycle on $Z$.  First consider $C\eta$ on
$CZ$, as a purely formal object.  Its main property is that
$\omegaCZ\frown C\eta$ is equivalent to $\setinfty\times\eta$.  In
addition, if $\eta$ and $\xi$ on $Z$ have complementary dimensions (and so
intersect to give a number), then $C\eta\frown\setinfty\times\xi =
\eta\frown\xi$.  These properties are not enough, in general, to determine
$C\eta$ as a homology class on $CZ$.  When $\xi$ is upper, $\omega$ is
injective, and so $\omega\frown\xi$ determines $\xi$.  This does not help
in the other dimensions.  Here, the Lefschetz isomorphism will be used. 
Write $\eta$ as $\omegaZ^r\eta'$, where $\eta$ and $\eta'$ have
complementary dimensions.  This representation is always possible and
unique, provided $\eta$ is lower.  Provided $r>0$, one can take
$\setinfty\times (\omegaZ^{r-1}\eta')$ as the cycle on $CZ$ that
represents the formal object $C\eta$.  Between them, these two cover all
the cases.  Thus, there is no formal obstacle to thinking of $C\eta$ as a
homology cycle on $CZ$.

The following construction, at least in certain cases, leads to a
geometric form for $C\eta$, in the lower dimensions.  As motivation, think
of $C\bfP_n$ (the cone on $\bfP_n$, not complex projective $n$-space). 
Here, each cycle $\eta$ on $\bfP_n$ can be coned to give a cycle $C\eta$
on $C\bfP_n$, \emph{provided the apex is not part of the stratification}. 
Adding the apex as a stratum will not however change the homology.  It
will thus be possible to move $\eta$ a little bit, so that it avoids the
apex (at least in the lower dimensions).  This change can be confined to a
small ball centered at the apex.  Think now of $C\eta$ on $CZ$ as follows.
Cone $\eta$ to form $C\eta$, and form a small ball $B$ about the apex. 
Outside of $B$, there is no fault with $C\eta$.  The task now is to change
$C\eta$ within $B$, so as to avoid the apex.  Consider now the
intersection $R=S\cap C\eta$ of $C\eta$ with the boundary $S$ of $B$.  For
certain values of $R\subset S$, it will be possible to `fill-in' $R$
within $B$, to obtain part of an intersection homology cycle, and for
others it will not.  (The difference between any two solutions is,
clearly, a local intersection homology cycle.)  For certain $\eta$ it will
be possible to solve the associated $R\subset S$ problem.  For heuristic
purposes, it will be assumed that this is always possible.  Intersection
properties can be used to resolve the indeterminacy due to local cycles at
the apex.  In this way it is possible (modulo some assumptions) to treat
$C\eta$ as a cycle on $CZ$, when $\eta$ is any mpih cycle on $Z$.  The key
is to if necessary modify the geometric form of $C\eta$ within a small
ball centered at the apex.

Now consider $C\eta$ on $CZ$, where $\eta$ is a local-global cycle.  The
example of $CCIC\sqdot$ forces one to allow this cycle, in some form of
the other.  The previous paragraph shows how this might be done.  The
geometric form of $C\eta$ must be modified in a small ball centered about
the apex of $CZ$, or perhaps more exactly, replaced by something else. 

In fact, this $C\eta$ problem need only be solved for mpih cycles.  Each
local-global cycle can be thought of as a $\lambda$-family of cycles.  One
can then solve this problem in the simpler case of $\lambda=1$, and then
define $C\eta$ to be the result of applying scalar multiplication to the
this solution.

The third type of cycle on $CZ$ are those obtained by moving a cycle
$\setinfty\times\eta$ away from the base $Z$, and then `coning it away'.
This was the point of departure, for the local-global theory.  (The mpih
theory prohibits the use of such objects, to generate homology relations. 
The local-global theory treats such objects as a cycle, but of a new
type.)  This construction can be iterated.  Here is an example. First, let
$\eta$ be the local-global cycle on $CIC\sqdot$.  Now let $Z$ be
$IICIC\sqdot$, and on $II\sqdot$ let $\xi$ be a coprimitive cycle.  On $Z$
there is a local-global cycle that can be written as $\xi\otimes\eta$. Now
consider $CZ$.  Provided $\setinfty\times(\xi\otimes\eta)$ can be moved
within its class, so as to avoid the base $Z$, it can be `coned away'. 
This is an example of a second-order local-global intersection homology
cycle.

Earlier in this section, $h\Delta$ was presented for all the $IC$
polytopes of dimension at most $4$.  To conclude, much the same will be
done for dimension $5$.  However, to save space this will be done
only for certain polytopes, whose flag vectors provide a basis for all
polytope flag vectors.  They are the ones in which neither $II$ nor
$I\sqdot$ occur.  There are $F_6=8$ such polytopes.

For these basis polytopes the $h$-vectors are as follows.
\[
\begin{array}{rl}
h(CCCCC\sqdot) &= (111111) \\
h(CCCIC\sqdot) &= (122221) +(111)\{1\} +(11)A\{1\} + (1)AA\{1\} \\
h(CCICC\sqdot) &= (122221) +\phantom{(111)\{1\}}
                                       +(11)A\{1\} + (1)AA\{1\} \\
h(CICCC\sqdot) &= (122221) +\phantom{(111)\{1\}}
                                +\phantom{(11)A\{1\}} + (1)AA\{1\} \\
h(CICIC\sqdot) &= (134431) +(111)\{1\}
                                +(11)A\{1\} + (2)AA\{1\} + (1)\{2\}\\
h(ICCCC\sqdot) &= (122221) \\
h(ICCIC\sqdot) &= (134431) +(121)\{1\} +(12)A\{1\} + (1)AA\{1\} \\
h(ICICC\sqdot) &= (134431) +\phantom{(111)\{1\}}
                                +(11)A\{1\} + (1)AA\{1\}
\end{array}
\]

Here is a summary of the discussion of the cycles on $IZ$ and $CZ$.  On
$IZ$ the cycles are as given by the K\"unneth principle.  Each cycle on
$IZ$ is a sum of products of a cycle on $I$ (or $\bfP_1$) with a cycle on
$Z$.  For $CZ$ the situation is more complicated.  If $\eta$ is a cycle on
$Z$, then one always has $\setinfty\times\eta$ on $CZ$.  Provided the
details are satisfied, as to what happens near to the apex, one will also
have $C\eta$.  (When $\eta$ is upper, these details are vacuous.) Finally,
if $\setinfty\times\eta$ can be moved so as to avoid the base $Z$ of $CZ$,
one also has its `coning away', the $\lambda$-cone $C_\lambda\eta$.  (A
necessary, and perhaps sufficient, condition for doing this is that
$\omegaCZ\frown\setinfty\times\eta$, which is equal to
$\setinfty\times(\omegaZ\frown\eta)$ be homologous to zero.)

This description of cycles motivates both the definition of the extended
$h$-vector (\S4), and the topological definition of local-global homology
(\S5).  In both cases, there are two aspects to the discussion.  The first
is the cycles themselves, the focus of this section.  The second is how
they are to be counted. 

Consider once again $ICIC\sqdot$ and $CCIC\sqdot$.  There, it was seen
that the same cycle may contribute to several different parts of
$h\Delta$.  This is something that is quite new.  The basic idea is this. 
A local-global homology group is spanned by all cycles that satisfy
certain conditions.  If these conditions are relaxed, another local-global
group is obtained.  (The same happens in intersection homology, when the
perversity is relaxed.) This is why the same local-global cycle may
contribute to several components of $h\Delta$.  The conditions are related
to where the cycle may be found. For example ($\dim=4$), one can count all
local $1$-cycles (subject to equivalence), or one can allow only those
that have some degree of freedom, as to their location.  The former are
counted by $\{1\}A$, the later by $x\{1\}$.  For example, $CICC\sqdot$ has
the former but not the later, while $ICIC\sqdot$ has both.  These
conditions control both the cycles and the relations.  Sometimes (the
$4$-cross polytope for example), relaxing the conditions may allow new
relations to appear amongst existing cycles. The varieties produced by $I$
and $C$ are special, in that this never happens.  This makes the
computation of their $h$-vectors much easier. This fact is exploited by
the next section.

\section{The extended $h$-vector}

This section defines, for every convex polytope $\Delta$, an extended
$h$-vector $h\Delta$.  It does this by using rules $\Itilde$ and $\Ctilde$
that satisfy the $IC$ equation.  These rules are motivated by the previous
section.  In the next section, local-global homology groups will be
defined for algebraic varieties.  Provided various assumptions are
satisfied, for $\Delta$ an $IC$ polytope the extended $h$-vector $h\Delta$
will give the local-global Betti numbers of the associated toric variety
$\PDelta$.  (The same may not be true for other rational polytopes, and
even it true will most likely be much harder to prove.  Such would be both
a formula for the local-global Betti numbers, and a system of linear
inequalities on the flag vectors of rational polytopes.)

There are two stages to the definition of $h\Delta$.  The previous section
described the local-global cycles on $IZ$ and $CZ$ in terms of those on
$Z$.  It also noted that the same cycle might contribute in several ways
to the local-global homology.  The first stage is to define operators
$\Itilde$ and $\Ctilde$ that count the local-global cycles, but without
regard to the multiple contributions.  The second stage is to make a
change of variable, to accomodate the multiple contributions.  This
corresponds to knowing the implications among the various conditions
satisfied by local-global cycles.  The second stage is vacuous for the mpih
paart of the theory.  (For the $IC$ polytopes, each
local-global cycle is determined by the corresponding global cycle,
together with a statement, as to the $\lambda$-coning conditions it
satisfies.  The first stage counts each local-global cycle only once, at
the most stringent conditions it satisfies.  This process is meaningful only
for $IC$ and similar polytopes.  In general, relaxation of conditions will
admit new relations, as well as new cycles.)

The first stage is to introduce an auxiliary vector $\htilde\Delta$,
defined via rules $\Itilde$ and $\Ctilde$.  The quantity $\htilde\Delta$
will be a sum of terms such as $[abcd]W$.  As noted in \S2, $[abcd]$
stands for the homogeneous polynomial $aX^3+bX^2Y+cXY^2+dY^3$, while $W$
will be a word in $\{k\}$ and $\Abar$.  The rules $\Itilde$ and $\Ctilde$
will be defined by their action on such terms.

The rule for $\Itilde$ is to multiply by $[11]=X+Y$.  It corresponds to
the K\"unneth formula for cycles.  The equation
\begin{equation}
\label{eqn.Itilde}
    \Itilde [abcd] W = [11] [abcd] W = [a,a+b,b+c,c+d,d] W
\end{equation}
is an example of this rule.  Note that this rule preserves the property of
being palindromic.  If swapping $X$ and $Y$ leaves $\htilde\Delta$
unchanged, then the same is true of $\Itildehtilde\Delta$.

The rule for $\Ctilde$ is more complicated.  It has three parts.  The
first part $\Ctilde_1$ leaves $W$ unchanged.  It corresponds to the idea,
that if $\eta$ on $Z$ is primitive, then so is $C\eta$ on $CZ$.  Here are
some examples of the rule
\[\begin{array}{ll}
    \Ctilde_1 [a]W = [aa]W \>;\quad &
    \Ctilde_1 [ab]W = [aab]W \>;\\
    \Ctilde_1 [abc]W = [abbc]W \>;\quad &
    \Ctilde_1 [abcd]W = [abbcd]W \>;\\
    \Ctilde_1 [abcde]W = [abccde]W \>;\quad &
    \Ctilde_1 [abcdef]W = [abccdef]W \>;
\end{array}\]
for this part.  It is to repeat the exactly middle, or failing that the
just before middle, term in the $[\ldots]$ sequence.

The reader is asked to verify that the equation
\[
    (\Itilde - \Ctilde_1) \Ctilde_1 = [010]
\]
holds, in the sense the applying the left hand side to, say, $[abcde]$
will produce $[0abcde0]$.  As $[11]$ and $[010]$ commute, $\Itilde$ and
$\Ctilde_1$ satisfy the $IC$ equation.  Using just this part of the rule
for $\Ctilde$ (together with the rule for $\Itilde$, and $h(\sqdot)=[1]$
as an initial condition) will generate the mpih part of $h\Delta$.

The second part $\Ctilde_2$ corresponds to the $\lambda$-coning of a cycle $\eta$ on
the base.  Such a cycle $\eta$ must be coprimitive.  The numbers $b'=b-a$,
$c'=c-b$ and so on count the coprimitives.  New words will be obtained by
prepending to $W$ a record, in the form $\Abar^j\{k\}$, of the coprimitive
that has been $\lambda$-coned.  Here are some examples of the rule
\[\begin{array}{ll}
    \Ctilde_2 [a]W = 0 \>;\quad &
    \Ctilde_2 [ab]W = 0 \>;\\
    \Ctilde_2 [abc]W = [b']\{1\}W \>;\quad &
    \Ctilde_2 [abcd]W = [b']\Abar\{1\}W \>;\\
    \Ctilde_2 [abcde]W = [b']\Abar^2\{1\}W +[c']\{2\}W 
                    \>;\quad \hidewidth \\
    \Ctilde_2 [abcdef]W = [b']\Abar^3\{1\}W +[c']\Abar\{2\}W 
                    \>;\quad \hidewidth 
\end{array}\]
for this part.  In $\Abar^j\{k\}$ the $k$ records the degree of the
coprimitive, and the $j$ `takes up the slack', to ensure homogeneity.  As
noted the previous section, the trivial coprimitives (which correspond to
$a'=a$) are not counted. 

The sum $\Ctilde_1 + \Ctilde_2$ of these two parts is not enough (or more
exactly, is too much).  One reason is that when used with $\Itilde$, the
result does not satisfy the $IC$ equation.  The third part is a
correction, that balances the books.  It is to subtract $[a]\Abar^j$, for
the appropriate power of $j$.  The geometric meaning of this correction
will be presented later.

Here now is the rule for $\Ctilde$. The examples
\begin{equation}
\begin{array}{rl}
    \Ctilde[a] &= [aa] - [a]\Abar                    \\
    \Ctilde[ab] &= [aab] - [a]\Abar^2                \\
    \Ctilde[abc] &= [abbc] - [a]\Abar^3  +[b']\{1\}              \\
    \Ctilde[abcd] &= [abbcd] - [a]\Abar^4  +[b']\Abar\{1\}       \\
    \Ctilde[abcde] &= [abccde]-[a]\Abar^5+[b']\Abar^2\{1\}+[c']\{2\}\\
    \Ctilde[abcdef] &= [abccdef] -[a]\Abar^6 +
                        [b']\Abar^3\{1\} +[c']\Abar\{2\} \\
\end{array}
\end{equation}
suffice to show the general rule.  In the above, it is to be understood
that both sides have been multiplied on the right by a word $W$ in the
symbols $\Abar$ and $\{k\}$.

This rule, and the rule for $\Itilde$, together satisfy the $IC$ equation. 
Here is an example.  The calculation
\[
    \Itilde\Ctilde [abcde] = [11][abccde] - [aa]\Abar^5
                    +[b'b']\Abar^2\{1\}
                    +[c'c']\{2\}
\]
follows immediately from the above.  The calculation for
$\Ctilde\Ctilde[abcde]$ is more involved.  One has
\[
    \Ctilde[abccde] = [abcccde] - [a]\Abar^6  
                    +[b']\Abar^3\{1\}
                    +[c']\Abar\{2\}    
\]
and also
\begin{eqnarray*}
    -\Ctilde[a]\Abar^5 &=& -[aa]\Abar^5 + [a]\Abar \Abar^5 \\
    \Ctilde[b']\Abar^2\{1\} &=& [b'b']\Abar^2\{1\} - [b']\Abar \Abar^2\{1\} \\
    \Ctilde[c']\Abar^2\{1\} &=& [c'c']\{1\} - [c']\Abar\{1\}
\end{eqnarray*}
as the various contributions.  Now compute the difference.  All but two of
the terms cancel.  One has
\[
    (\Itilde\Ctilde-\Ctilde\Ctilde)[abcde] = [11][abccde] - [abcccde]
\]
which is, as for $\Ctilde_2$, is equal to $[010][abcde]$.  As this example
is completely typical, the result follows.

To complete this definition of $\htilde\Delta$, one must supply an initial
value $\htilde(\sqdot)$, such that $\Ctilde\htilde(\sqdot)$ and
$\Itildehtilde(\sqdot)$ are equal.  Here a problem arises.  The value
$\htilde(\sqdot)=[1]$ does not quite work.  The quantities
$\Ctilde[1]=[11]-[1]\Abar$ and $\Itilde[1]=[11]$ are not equal.  Here is
the solution.  Recall that $h\Delta$ is to be a sum of terms of the form
$h_W\Delta\cdot W$, for some family of symbols~$W$. Thus, one should
really write $\htilde(\sqdot)=[1]W_0$, for some symbol $W_0$. Consistency
is then equivalent to the equation $\Abar W_0 =0$.  However, it is more
convenient to use `$\sqdot$' as the initial value.  The initial values
\begin{equation}
    \htilde (\sqdot) = [1] \sqdot \> ; \qquad \Abar \sqdot = 0
\end{equation}
conclude the definition.  The multiple use of the symbol `$\sqdot$' in
practice causes no confusion.

(The geometric meaning of the correction $-[a]\Abar^i$ is as follows. 
Suppose that $\Delta$ is an $IC$ polytope, and that say $[abc]W$ appears
in $\htilde\Delta$.  This term $[abc]W$ is due to the $I$ and $C$
operations that constructed $\Delta$.  In particular, going back in this
case two steps, one obtains the polytope $\Delta_1$, from which $\Delta$
is derived by applying $I$ and $C$, and in $\htilde\Delta_1$ the term
$[a]W$ appears.  Going back another step, one has $\Delta_1=C\Delta_0$,
and on $\Delta_0$ there are $a$ independent coprimitive cycles $\eta$,
that contribute $[a]W$ to $\htilde\Delta_1$.  Now consider $C\Delta$.
These coprimitive cycles $\eta$ (multiplied by $\setinfty$ as appropriate)
continue to exist on $\Delta$, and the $\Ctilde_2$ part of the rule for
$\Ctilde$ will thus cause them to contribute $[a]\Abar^iW$. However, their
contribution to $\htilde C\Delta$ has already been counted, as the `$a$'
part of $\Ctilde_1[abc]W=[abbc]W$.  This is because $\htilde\Delta$ counts
cycles on $IC$ polytopes, at the most stringent condition they satisfy.
Hence the correction $-[a]\Abar^i$.)

Finally, note that if $\htilde\Delta$ is palindromic then so is
$\Ctilde\htilde\Delta$.  As the same is true for $\Itilde$, it follows
that $\htilde\Delta$ is indeed palindromic, first for all $IC$ polytopes
and then (by linearity) for all polytopes.  This completes the definition
of the auxiliary vector $\htilde\Delta$.

The second stage in the definition of the extended $h$-vector is to apply
a linear change of variable to the above quantity $\htilde\Delta$.  The
meaning of this transformation is as follows.  Recall that for $\Delta$ an
$IC$ polytope, $\htilde\Delta$ counts each local-global cycle once,
according to the most stringent conditions it satisfies.  The
transformation is the process of relaxation of conditions, as it applies
to cycles on $IC$ polytopes.

Previously, little attention has been given to the conditions satisfied by
the local-global cycles.  These cycles have been generated from
coprimitives via the $\lambda$-coning operation, but the topologicial
properties of the cycles so obtained have not been explicitly formulated. 
Instead, the $X$, $Y$, $\{k\}$ and $\Abar$ symbols have been used in a
somewhat formal way, to record the relevant facts relating to the
construction of the associated cycles.

Recall that in the discussion of mpih, one first obtained the cycles, and
then the conditions that they satisfied.  Something similar will be done
with local-global homology.  A local-global cycle on an $IC$ variety
$\Delta$ consists of a cycle, as counted by $\htilde\Delta$, together with
a condition that it satisfies.  Such conditions, which control how the
apex locus meets the strata, will be formulated in the next section. 
However, by assumption each word in $x$, $y$, $A$ and $\{k\}$ determines a
condition, or in other words a type of local-global cycle.  The present
task is to describe how to pass from $\htilde\Delta$ to $h\Delta$.  If
done properly, it will implicitly determine the set of conditions, that
will in the next section be explicitly stated.

To each term such as $\Abar\{1\}$ or $X\{1\}$ in $\htilde\Delta$, there is
to correspond a condition that applies to local-global cycles.  The
coefficient in $\htilde\Delta$ of such a term counts how many cycles there
are on $\Delta$ that satisfy that condition, but not any condition that is
more stringent.   (This makes sense only for $IC$ polytopes.  Other
polytopes may produce negative coefficients.) The transformation from
$\htilde$ to $h$ is thus determined once one knows the partial order on
conditions, that is associated to stringency or implication.  In other
words, the partial order is that of implication among the associated
conditions.

Here is an example.  In $\htilde(CCIC\sqdot)$ the term $X\{1\}$ occurs,
with coefficient~$1$.  (As noted in \S3, it also so occurs in
$\htilde(ICIC\sqdot)$.) This term corresponds to a particular type of
local-global cycle, namely one that can be found anywhere along the apex
edge.  Now consider $\htilde(CICC\sqdot)$.  Here the term $\Abar\{1\}$
occurs, with coefficient~$1$.  This too corresponds to a local-global
cycle, but of a different type. It corresponds to a local-global cycle
that can be found only about the apex of $CICC\sqdot$.  (In the previous
case, there was an apex edge, anywhere along which the local-global cycle
could be found.  In this case there is no apex edge.)  Thus the condition
$X\{1\}$ is more stringent than $\Abar\{1\}$.

Each term such as $\Abar\{1\}$ or $X\{1\}$ in $\htilde\Delta$ will
contribute $A\{1\}$ and $x\{1\}$ respectively to $h\Delta$.  These
`diagonal' terms arise because `the most stringent condition satisfied by
a cycle' is also `a condition satisfied by a cycle'.  In addition, as just
noted, if a cycle satisfies $X\{1\}$ then it will also satisfy
$\Abar\{1\}$ (but not as the most stringent condition), and so $X\{1\}$
will contribute $A\{1\}$, in addition to $x\{1\}$.  The transformation of
\[
    \htilde (CCIC\sqdot) = [12221] + [11]\{1\}
\]
is therefore
\[
          h (CCIC\sqdot) = (12221) + (1)A\{1\} + (11)\{1\} 
\]
which is as postulated in the previous section.

Here is another example.  In involves second order local-global cycles. 
The symbol $\{1\}\{1\}$ represents a local $1$-dimensional family of local
$1$-cycles.  In dimension~$7$ each of the terms $\Abar\{1\}$,
$X\{1\}\{1\}$ and $\{1\}\Abar\{1\}$ is to represent a different condition
on a local $1$-family of local $1$-cycles.  To begin to understand these
conditions, for each of these terms an $IC$ polytope will be produced, in
whose $\htilde$-vector the term occurs.  Here is a list
\[\begin{tabular}{l}
    in $\htilde(CICCCIC\sqdot)$ the term $\Abar\{1\}\{1\}$ occurs \\
    in $\htilde(CCICCIC\sqdot)$ the term $X\{1\}\{1\}$ occurs \\
    in $\htilde(CICCICC\sqdot)$ the term $\{1\}\Abar\{1\}$ occurs 
\end{tabular}\]
of such polytopes.  This the reader is invited to verify.

From this list the conditions, or more exactly the partial order, will be
obtained.  As $IC$ varieties, each of the above will have a minimal
stratification.  The conditions control how the apex loci of the cycle
meets the strata.  (The apex locus is due to the $\lambda$-coning or
`locality' of the cycle.) Here is a list, which again the reader is
invited to verify,
\[\begin{tabular}{l}
    the strata for $\Abar\{1\}\{1\}$ have dimension $0$, $4$ and $7$ \\
    the strata for $X\{1\}\{1\}$ have dimension $1$, $4$ and $7$ \\
    the strata for $\{1\}\Abar\{1\}$ have dimension $0$, $3$ and $7$ 
\end{tabular}\]
of the strata dimensions for the associated $IC$ varieties.

The present discussion is of a local $1$-dimensional family of local
$1$-cycles, and so there are two apex loci to consider.  One of them lies
on the strata of dimension $0$ or $1$, the other on strata of dimension
$3$ or $4$.  It is more stringent to require that an apex locus be found
on a $1$-strata than on a $0$-strata, and similarly for $4$-strata and
$3$-strata.  Thus
\[
    X\{1\}\{1\} \implies \Abar\{1\}\{1\} \implies \{1\}\Abar\{1\}
\]
is the partial order on conditions.  As in the example of $\Abar\{1\}$ and
$X\{1\}$, this partial order allows $\htilde\Delta$ to be transformed into
$h\Delta$.

It is now possible to define the partial order on terms.  As in the
previous example, associate to each term an $IC$ polytope, whose
$\htilde$-vector realises the term.  Each such polytope, thought of as a
variety, has a minimal stratification.  The partial order on terms is then
the partial order on the dimension vectors of the stratification
associated to each term.

Here are the details.  The partial order does not compare, say,
$\{1\}\{2\}$ and $\{2\}\{1\}$, or $X\{1\}$ and $Y\{1\}$.  The associated
cycles differ by more than a change in conditions.  Say that two
expressions in $X$, $Y$, $\Abar$ and $\{k\}$ are \emph{broadly similar}
if, when $X$ and $\Abar$ have been deleted, they are identical.  They
should also have the same degree.  The partial order applies only to
broadly similar terms.  Each term will appear in $\htilde\Delta$, for one
or more $IC$ polytopes $\Delta$. Choose one of these polytopes (it does
not matter which) for which the associated stratification has as few terms
as possible.  It will in fact have $r+1$ terms, where $r$ is the order of
the term.  This process associates a stratification dimension vector
$d=(d_1<d_2\ldots<d_{r+1}=n)$ to each term.  The partial order is that $d$
will imply $d'$ just in case $r=r'$, and $d_i\geq d_i'$ for $i=1,2,\ldots,
r+1$.

The above is a description, based on geometry, of the partial order on
conditions.  The stratification dimension vector $d$ associated to a term
can also be computed directly.  Suppose, for example, the term is
$X^2Y^3\Abar^4\{5\}\Abar^2\{6\}$.  The degree of this term is
\[
    2 + 3 + 4 + (2\times5+1) + 2 + (2\times6+1) = 35
\]
and so $35$ is the dimension of the top stratum.  For the broadly similar
term $X^i\{5\}\{6\}$, where $i$ makes the degree up to $35$, the dimension
vector written backwards is
\[
    35 > 35 - \deg \{6\} > 35 - \deg \{6\} - \deg \{5\}
\]
(or $11 < 22 < 35$ the normal way around).  For the given term
$X^2Y^3\Abar^4\{5\}\Abar^2\{6\}$ the dimension vector is 
\[
    35 > 35 - \deg \{6\} - 2 > 35 - \deg \{6\} - 2 - \deg \{5\} -4
\]
and from this the partial order can be given a combinatorial description. 
(The proofs have been left to the reader.)

The partial order on terms can be expressed in the following way.  The
terms implied by some given term, say $X^2Y^3\Abar^4\{5\}\Abar^2\{6\}$,
can all be obtained in the following way.  First of all, any number of
occurences of $X$ can be replaced by $\Abar$.  Such $\Abar$ should of
course be placed after the $X^iY^j$ term.  Next, any occurence of $\Abar$
can be `slid' rightwards over any $\{k\}$ symbol.  Finally, note that the
terminating symbol `$\sqdot$' is assumed to be at the end of each word,
and that $\Abar\sqdot$ is zero.

It is now possible to give an algebraic description of the rules that
transform $\htilde\Delta$ to the extended $h$-vector $h\Delta$.  First of
all the change of variable
\begin{equation}
\label{eqn.XxYy}
    X = x + \Abar \>;\quad Y=y \>;\quad
\end{equation}
is made.  To ensure that equations such as
\[\begin{array}{l}
    X^2 = x^2 + x\Abar + \Abar^2 \> ; \\
    X^3 = x^3 + x^2\Abar + x\Abar^2 + \Abar^3 \> ;
\end{array}\]
hold, the equation
\begin{equation}
\label{eqn.Abar-x}
    \Abar x = 0 
\end{equation}
is postulated.  These rules allow $X$ and $Y$ to be replaced by $x$ and
$y$.  Of course, $X$ and $Y$ commute, as do $x$ and $y$.  Thus, $y$ and
$\Abar$ are also to commute.

The equation
\begin{equation}
\label{eqn.Abar-i}
    \Abar \{k\} = A\{k\} + \{k\}\Abar
\end{equation}
allows the `rightward slide' of $\Abar$ over $\{k\}$.  The $\Abar$ on the
left hand side can remain where it is, to produce $A\{k\}$; or it can
slide over the $\{k\}$.  If so slid, it could be slid again.  Thus, it is
$\{k\}\Abar$ rather than $\{k\}A$ on the right hand side.  Finally, in
$\Abar A$ the second $A$ is `non-sliding'.  This stops the $\Abar$ from
sliding.  Thus the equation
\begin{equation}
\label{eqn.Abar-A}
    \Abar A = AA
\end{equation}
is also postulated.

The rules in the previous two paragraphs define the extended $h$-vector
$h\Delta$ of any convex polytope $\Delta$, via the auxiliary vector
$\htilde\Delta$.  Implicit in this are topological conditions on
local-global cycles.  In the next section these conditions will be made
explicit.

\section{Topology and local-global homology}

In this section $Z$ will be a complex algebraic variety, considered as a
stratified topological space.  Local-global intersection homology groups
will be defined for $Z$.  The starting point is a particular means of
expressing the concept of a local cycle, and of course the basic concepts
of intersection homology.  By considering families of such cycles, the
full concept of a local-global cycle is developed.  This is done by
extending the notion of a simplex.

Recall that in the previous section there were two stages to the
definition of the $h$-vector $h\Delta$.  Loosely speaking the first stage,
the definition of $\htilde\Delta$, corresponded to the definition of a
local-global cycle.  The second stage was concerned with the conditions
(on how it meets the strata) satisfied by the cycle.  The topological
definition of local-global homology similarly has two stages, the
definition of the cycles, and the definition of the conditions.  Most of
the justification for the cycle definition comes from \S3, while the
conditions are those implicit in \S4.

To begin with, consider the concept of a local cycle.  In \S2 such was
thought of as a cycle $\eta$ lying on any sufficiently small neighbourhood
$U$ of the point $s$, about which the cycle is to be local. That the local
topology of $Z$ about $s$ is the cone on the link can clarify this
concept.  In this section, all this will be incorporated into the
definition of a local cycle.  Some preliminaries are required.

Homology can be defined using embedded simplices.  Here is a review.  An
\emph{embedded $k$-simplex} is simply a continuous map $f:\sigma_k\to Z$
from the $k$-simplex $\sigma_k$.  A \emph{$k$-chain} is a formal sum of
(embedded) $k$-simplices.  Each simplex $\sigma_k$ has a boundary
$d\sigma_k$, which is a formal sum of $(k-1)$-simplices.  In the same way,
each $k$-chain $\eta$ has a \emph{boundary} $d\eta$, which is a
$(k-1)$-chain.  A \emph{$k$-cycle} is a $k$-chain $\eta$ whose boundary
$d\eta$ is zero (as a formal sum of embedded simplices).  Because $d\circ
d = 0$, if $\xi$ is a chain then $\eta=d\xi$ is a cycle.  Such a cycle
$\eta$ will be called a \emph{boundary}.  The \emph{$k$-th ordinary
homology group} of $Z$ consists of the $k$-cycles modulo the
$k$-boundaries.

Intersection homology \cite{bib.MG-RDM.IH} places restrictions on how the
embedded simplices can meet the strata.  The conditions are expressed
using a perversity, which is a sequence of numbers.  In this paper only
the middle perversity is used.  It imposes the following condition on an
embedded $k$-simplex $f:\sigma_k\to Z$.  Let $S_i$ be the complex
$i$-dimensional stratum of $Z$.  Consider $f^{-1}(S_i)$, or more exactly
its dimension.  Let the empty set have dimension $-\infty$.  If the
inequality
\[
    \dim f^{-1}(S_i) \leq k - ( n-i ) 
\]
holds for every stratum $S_i$, then $f:\sigma_k\to Z$ is \emph{allowed}
(for the middle perversity).  An \emph{allowed cycle} $\eta$ is a formal
sum $\eta$ of embedded simplices, whose boundary is zero.  An
\emph{allowed boundary} $\eta=d\xi$ is a formal sum $\xi$ of allowed
simplices, whose boundary $d\xi$ (which necessarily is a cycle) is also
allowed.  The $k$-th (middle perversity) \emph{intersection homology}
group $z$ consists of the allowed cycles modulo the allowed boundaries. An
important technical result \cite{bib.MG-RDM.IH2,bib.HK.TIIH} is that this
group is independent of the stratification $S_i$ chosen for $Z$.

Consider the concept of a local cycle.  Each point $s$ on $Z$ has a
neighbourhood that is homeomorphic to the cone $CL_s$ on something, namely
the \emph{link} $L_s$ at $s$.  One can use $CL_s$ as the neighbourhood $U$
in which the cycles $\eta$ local to $s$ can be found. Because $CL_s$ is a
cone (and $\eta$ avoids the apex of the cone), the cycle $\eta$ is
equivalent to a cycle $\eta'$, that is supported on the base $L_s$ of the
cone.  (Use the cone structure to move the cycle $\eta$ to the base of the
cone.) In fact, for each $0<\lambda\leq 1$ one obtains a cycle
$\eta_\lambda$, while the `limit' $\eta_0$ of the family is the apex
$CL_s$, which is the point $s$ of $Z$.  In other words, a local cycle is a
cycle that can be `coned away' to a point, except that the perversity
conditions may disallow this.

The concept of a local cycle can be formulated without using either $U$ or
$CL_s$.  Define a \emph{coned $k$-simplex} to be the cone $C\sigma_k$ on a
$k$-simplex.  (Although isomorphic to $\sigma_{k+1}$, it will not be
treated as such.  Later, more complicated objects will be coned.)  Local
homology will be constructed using such simplices.  An \emph{embedded
coned $k$-simplex} is simply a continuous map $f:C\sigma_k\to Z$.  The
\emph{apex} of $f$ is the image under $f$ of the apex of $C\sigma_k$.  A
\emph{coned $k$-chain} is a formal sum of (embedded) coned $k$-simplices. 
Each coned simplex $C\sigma_k$ has a boundary, which is a formal sum of
coned $(k-1)$-simplices.  (The boundary is taken only the the $\sigma_k$
direction, and not in the $C$ direction.)  In the same way, each coned $k$-chain $\xi$ has
a \emph{boundary} $d\xi$, which is a coned $(k-1)$-chain.  A \emph{coned
$k$-cycle} is a coned $k$-chain whose boundary $d\eta$ is zero (as a
formal sum of embedded coned $(k-1)$-simplices.)  Because $d\circ d = 0$,
if $\xi$ is a coned chain, then $d\xi$ is a coned cycle.  Such a cycle
$\eta=d\xi$ will be called a \emph{coned boundary}.

Finally, only certain cycles and boundaries will be allowed.  Use
$\lambda$ and $p$ to denote the cone and simplex variable respectively. 
For each $0<\lambda\leq 1$, use $f_\lambda$ to denote the embedded simplex
defined by the rule $f_\lambda(p)=f(\lambda,p)$, where $f$ is an embedded
coned simplex.  Say that a coned cycle $\eta$ is \emph{allowed} if
$\eta_\lambda$ is similarly allowed.  Say that a coned boundary
$\eta=d\xi$ is \emph{allowed} if $\eta_\lambda=d\xi_\lambda$ is allowed,
for $0<\lambda\leq1$.  Now fix a point $z\in Z$.  Say that an embedded
coned simplex $f$ is \emph{local to $s$} if $s$ is the apex of $f$.  In
that case, $s$ will also be the apex of the boundary of $f$. The
\emph{local $k$-homology at $s$} consists of the coned $k$-cycles modulo
the coned $k$-boundaries, where the cycles and boundaries are allowed (by
the perversity), \emph{and are constructed using only the embedded coned
simplices local to $s$}.

Local-global cycles differ from local cycles in that the base point or
apex is allowed to move.  To do this a new sort of simplex is required. 
Consider $\sigma_1\times C\sigma_k$.  This can be thought of as a
$1$-dimensional family of coned $k$-simplices.  It has not a single point
as its apex, but an \emph{apex edge}.  Its boundary in the $\sigma_1$
direction is a pair of coned $k$-simplices, one at each end of the apex
edge.  It also has a boundary in the $\sigma_k$ direction, which is a
formal sum of $\sigma_1\times C\sigma_{k-1}$ `simplices'. (As before, no
boundary is taken in the cone direction.)

The task now is to define the \emph{local-global $k$-homology} of $Z$. 
The cycles are as in local $k$-homology, except that it is not required
that the embedded coned $k$-simplices be local to some fixed point $s$ of
$Z$.  Clearly, such a cycle can be written as a formal sum of local
$k$-cycles, based at different points $s_1$, $\ldots$, $s_N$ of $Z$.

The boundaries require more thought.  Previously the simplices (possibly
coned) were indexed by a single number $k$, and the boundary operator $d$
reduced the index by one.  The present situation is that $C\sigma_k$
arises from both $C\sigma_{k+1}$ and $\sigma_1\times C\sigma_k$, and the
latter also produces $\sigma_1\times C\sigma_{k-1}$.  When $C\sigma_k$ is
written as $\sigma_0\times C\sigma_k$, it becomes clear that the
`simplices' are now indexed by a pair of numbers, and that the boundary
operator $d$ can reduce either one or the other by one.

Because of this, the concept of $\eta=d\xi$ being a \emph{local-global
$k$-boundary} can be formulated in several, possibly inequivalent, ways. 
In each case $\xi$ will be a formal sum of embedded $\sigma_i\times
C\sigma_j$ `simplices', with $i+j=k+1$.  One also requires that $\eta$ is
the boundary of $\xi$, and that for each $0<\lambda\leq 1$, one has the
$\eta_\lambda=d\xi_\lambda$ is allowed (in the usual way).  The different
concepts arise, according to the values of $i$ and $j$ allowed.  Say that
$\xi$ is \emph{very pure} if all its `simplices' have the same type. 
Clearly, very pure boundaries should be allowed.  In the present case,
$\sigma_0\times C\sigma_{k+1}$ is required for local equivalence, while
$\sigma_1\times C\sigma_k$ allows the apex to move.  The sum of two very
pure boundaries need not be very pure.  Say that $\xi$ is \emph{pure} if
all the `simplices' are capable, by virtue of their `dimension', of
participating in a very pure boundary.  In the present case, it means that
$\xi$ is built out of a mixture of $\sigma_0\times C\sigma_{k+1}$ and
$\sigma_1\times C\sigma_k$.  Finally, say that $\xi$ is \emph{mixed} if it
is built out of any mixture whatsoever of $\sigma_i\times C\sigma_j$
`simplices'.

In principle, the three concepts are different.  The author expects very
pure and pure to give the same boundaries.  (This means that if
$\eta=d\xi$ is a pure boundary, then there are very pure $\eta_i=d\xi_i$,
with $\eta=\sum\eta_i$.  It is not required that $\xi=\sum\xi_i$.)  The
mixed concept permits $\sigma_{k+1}\times C\sigma_0$ to play a r\^ole. 
This seems to be wrong.  In this paper, local-global homology will be
defined using the pure concept.  Experience will show if this is correct. 
(Already, the desired Betti numbers are known.)

The general definition of local-global cycles and boundaries proceeds in
the same way.  (Recall that the conditions that control how the apex loci
meets the strata have not yet been considered.)  To begin with, let
$k=(k_0,k_1,\ldots,k_r)$ be a sequence of nonnegative integers.  Define
the \emph{$k$-simplex} $\sigma_k$ to be the convex polytope
\[
    \sigma_k = \sigma_r \times C ( \sigma_{r-1} \times C ( \> \ldots\>
                (\sigma_1 \times C \sigma_0) \ldots ))
\]
where each $\sigma_i$ is an ordinary simplex of dimension $k_i$.  In other
words, $\sigma_k$ is $\sigma_r\times C\sigma_{k'}$, where
$k'=(k_0,k_1,\ldots, k_{r-1})$.  The number $r$ is the \emph{order} of
$\sigma_k$ (and of $k$).  It is also the number of coning operators.  If
$r=0$ then $k=(k_0)$, and $\sigma_k$ is an ordinary simplex, of dimension
$k_0$.

For each choice $\lambda=(\lambda_1, \ldots , \lambda_r)$ of a nonzero
value for each of the coning variables in $\sigma_k$, one obtains a
`simplex' $\sigma_{k,\lambda} \cong \sigma_r\times \ldots \times \sigma_0$.
A formal sum $\eta$ of continuous maps $f:\sigma_k\to Z$ is a
\emph{$k$-cycle} if the boundary $d\eta$ is zero, and for each such
$\lambda$ the maps $f_\lambda:\sigma_{k,\lambda}\to Z$ associated to
$\eta$ are allowed (by the perversity conditions).  The $k$-cycle $\eta$
is a \emph{$k$-boundary} if there is a formal sum $\xi$ of $k'$-simplices
such that $\eta=d\xi$, and again the $f_\lambda$ due to $\xi$ are also
allowed (by the perversity conditions). Here, because pure boundaries are
being used, the index $k'$ ranges over the $r+1$ indices obtained via
choosing one of the $k_i$ in $k$, and raising it by~$1$.  Note that the
boundary components of $\xi$, whose index is not $k$, must all cancel to
zero.

This defines the local-global cycles and boundaries.  Each local-global
homology group is determined by the choice of an index $k$ (which gives
the `dimension' of the cycles), and a choice of the conditions that
control how the cycles and boundaries meet the strata.  The rest of this
section is devoted to the study of these conditions.

Recall that $Z$ is a complex algebraic variety considered, as a stratified
topological space.  Each stratum $S_i$ has real dimension $2i$.  The basic
building block for local-global homology, the $k$-simplex $\sigma_k$, can
also be thought of as a stratified object.  However, instead of strata it
has \emph{apex loci}.  Altogether $\sigma_k$ will have $r$ (its order) apex
loci. The $k$-simplex $\sigma_k$ is produced using $r$ coning
operators. Each coning introduces an apex.  Multiplication by the
$\sigma_i$, and the subsequence coning operations, will convert this apex
into an apex locus; except that any new apex does not belong to the
already existing apex locus.  The closure of the $i$-th apex locus is a
$\sigma_{k'}$-simplex, where $k'=(k_i,k_{i+1}, \ldots, k_r)$.

The conditions, that control how the apex loci meet the strata, are to
give rise to groups that are independent of the stratification.  From
this, the nature of the conditions can be deduced.  Here is an example. 
Suppose $A$ and $B$ are nonsingular projective varieties.  Let $Z=A\times
CB$ be the product of $A$ with the (projective) cone on $B$, with $Z$ and
the apex locus $A$ as the closures of the strata.  Each coprimitive $\eta$
on $B$ will (if it can be moved away from the base) determine a
local-global cycle $C_\lambda\eta$ on $CB$ and then, by K\"unneth, the
choice of a cycle $\xi$ on $A$ determines a local-global cycle $\xi\otimes
C_\lambda\eta$ on $Z$.  Suppose now that $\xi$ is the fundamental class
$[A]$ on $A$, and that some condition permits $[A]\otimes C_\lambda\eta$. 
The apex locus of this cycle is the apex locus $A$ of $Z$.  Now add strata
to $Z$, that lie inside the apex locus $A$.  There is no possibility of
moving $[A]\otimes C_\lambda\eta$ within its equivalence class, in such a
way that the apex locus of the cycle is changed.  Thus, this addition of
new strata to $A$ will not affect the admissability of $[A]\otimes
C_\lambda\eta$.

Return now to the general case.  The previous example can be expressed in
the following way.  Use the stratification $S_i$ to $Z$ to define the
filtration
\[
    U_i = Z - S_0 - S_1 \ldots -S_{i-1} 
        = S_i \cup S_{i+1} \cup \ldots \cup S_n
\]
of $Z$ by open sets $Z=U_0 \supseteq U_1 \ldots \supseteq U_n$.  For each
apex locus $A$ of a local-global cycle, ask the following question: what is
the largest $i$ such that $A\cap U_i$ is dense in $A$?  Call this the
\emph{$w$-codimension} $w(A)$ of $A$.  (It tells one \emph{where} $A$ is
generically to be found.)  Clearly, adding strata as in the previous
example will not change $w(A)$.

Consider once again $Z=A\times CB$.  Suppose one wants a local-global
cycle on $Z$, whose $w$-codimension has some fixed value $i\leq\dim A$. 
Such can be achieved, provided $Z$ is suitably stratified.  To begin with,
let $A_i\subseteq A$ be a subvariety of dimension $i$.  Now form the
local-global cycle $[A_i]\otimes C_\lambda\eta$, where $[A_i]$ is the class
of $A_i$ in $H_i(A)$.  If $Z$ and $A$ are the only closures of strata,
then $w(A_i)$ will be $\dim A$.  However, one can always add $A_i$ to the
stratification, and then $w(A_i)$ will be equal to $\dim A_i$, which by
construction is equal to~$i$.

The meaning of this example is as follows.  Strata impose conditions on
cycles.  The more strata, the more conditions.  Suppose a condition allows
the $w$-codimension of the apex locus of a cycle to be some number $i$,
smaller than $\dim A$.  Cycles meeting this condition can be found, when a
suitable stratum (the subvariety $A_i$) is added.  Thus, the condition
must allow the cycle $[A_i]\otimes C_\lambda\eta$, even when $Z$ is given
its minimal stratification, if the result is to be stratification
independent.  Reducing the stratification can only reduce the
$w$-codimension.  Thus, if a condition allows $i$ as a $w$-codimension, it
should also allow all values smaller than $i$.  In other words, $i$ should
be maximum allowed value for the $w$-codimension.

(Another argument in favour of this conclusion is as follows.  The apex
locus of a local-global cycle need not be connected.  Suppose the apex
loci $A$ and $A'$ respectively of $\eta$ and $\eta'$ have $w$-codimensions
$i$ and $i'$, with $i<i'$.  Suppose also that $i$ is a permitted
$w$-codimension.  Thus, $\eta$ is permitted.  Now consider $\eta+\eta'$. 
Although its apex locus need not in every case be $A\cup A'$, so far as
$w$-codimension is concerned, it might as well be.  Thus, $\eta +\eta'$ will
also have $i$ as the $w$-codimension of its apex locus, and so is permitted. 
Because both $\eta$ and $\eta+\eta'$ are permitted, the difference
$(\eta+\eta')-\eta=\eta'$ must also be allowed.  In other words, $i$ is a
minimum allowed value.)

The \emph{local-global intersection homology groups} $H_{k,w}(Z)$ can now
be defined.  The subscript
\[
    k=(k_0,k_1,\ldots,k_r)
\]
is a sequence of positive
integers.  The cycles and boundaries of $H_{k,w}(Z)$ are constructed out
of embedded $k$-simplices, as defined earlier in this section.  The
\emph{$w$-condition} $w$ is a sequence 
\[
    w_1 \geq w_2 \geq \ldots \geq w_r \geq 0
\]
of nonegative integers, which is used as follows.  Each cycle (or
boundary) will have $r$ apex loci, to be denoted by $A_1$, $A_2$,
$\ldots$, $A_r$.  Each apex locus $A_i$ will have a $w$-codimension
$w(A_i)$.  Only those cycles and boundaries for which
\[
    w_i(A) \leq w_i \qquad i = 1, \ldots, r
\]
holds are to be used in the construction of $H_{k,w}(Z)$.  

(In this construction, it is assumed that the middle perversity is used. 
If the $\eta_\lambda$ and $d\xi_\lambda=\eta_\lambda$ are to satisfy some
other condition $p$, it can be added to the definition, and also to the
notation $H^p_{k,w}(Z)$, like so.  Only for certain values of $k$, $w$ and
$p$ will nonzero groups be possible.  The notation of the previous section
takes this into account.)

\section{Summary and Conclusions}

To close this paper, its main points will be summarized, and then various
questions arising are discussed.  These are firstly, are the local-global
intersection homology groups $H_{k,w}Z$ independent of the stratification
of $Z$?  Second, does $h\Delta$ compute the (local-global) Betti numbers
of $\PDelta$.  Third, can the (local-global) homology $H_\bullet\Delta$ of
$\Delta$ be constructed without recourse to $\PDelta$?  Fourth, what
happens when integer coefficients are used?  (This question is related to
the resolution of singularities.)  Fifth, does local-global homology have
ring- or functor-like properties.  Finally, there is a brief history of
the genesis of this paper, and acknowledgements.

The relative importance of the various parts of this paper depend on one's
point of view.  The construction (\S3) of local-global cycles on the $IC$
varieties was chosen as the starting point, from which both the formula
for $h\Delta$ (\S4) and the topological definition (\S5) followed.  For
one interested in more general varieties the topological definition is
perhaps most important.  The examples in \S3 then become merely the
application of a more general concept.  Finally, not all general polytopes
can be studied via topology, and so this gives $h\Delta$ (\S4) and the
questions arising from it a special importance.  In some sense,
$\lambda$-coning and apex loci are the key new concepts, over and above
perversity conditions on cycles and boundaries.  At the risk of pleasing
none, this paper has tried to please all.

The concept of local-global homology has two aspects, namely the cycles
and boundaries on the one side, and the strata conditions on the other. In
the notation $H_{k,w}Z$, it is $k$ that controls the type of cycles and
boundaries, while $w$ supplies the strata conditions.  A local-global
cycle $\eta$ is a global cycle $\eta_{(1)}$, which can be coned in any of
$r$ (the order of $k$ and of $\eta$) $\lambda$-directions.  The subscript
in $\eta_{(1)}$ indicates that $(1,1,\ldots,1)=(1)$ are the values
$\lambda_i$ of the coning variables.  The cycle $\eta$ is to be a suitable
family $\eta_\lambda$ of global cycles, which collapses in various ways as
each coning variable $\lambda_i$ goes to zero.

These $\lambda_i$ are not independent.  If $\lambda_j=0$ then
$\eta_\lambda$ does not depend on $\lambda_i$, for $i<j$.  In other words,
there is a sequence of \emph{collapsings}
\[
    \eta_0=\eta \gg \eta_1 \gg \eta_2 \gg \ldots \gg \eta_r
\]
of the family $\eta=\eta_\lambda$ of cycles.  Each $\eta_i$ is the result
of placing $\lambda_i=0$, and is called an \emph{apex locus}.  The
sequence $k=(k_0,k_1,\ldots,k_r)$ encodes the dimensions (more exactly the
relative dimensions) of these families.  These collapsings provide local
information about the cycle $\eta$.

Of interest is not only \emph{how} $\eta$ collapses to $\eta_i$, but also
\emph{where}.  By this is meant how $\eta_i$ meets the strata $S_j$ of
$Z$.  More exactly, $w(\eta_i)$ is defined to be the largest $j$ such that
\[
    \eta_i \cap ( S_j \cup S_{j+1} \cup \ldots \cup S_n)
\]
is dense in $\eta_i$.  This is a measure of how special a locus is
required, to generically effect the collapsing of $\eta$ to $\eta_i$.  The
index $w=(w_1,w_2,\ldots,w_r)$ in $H_{k,w}Z$ controls the cycles and
boundaries used as follows.  A cycle $\eta$ (or boundary $\xi$) is to be
used only if $w(\eta_i)$ (or $w(\xi_i)$) is at least $w_i$.  Thus,
reducing $w$ will potentially allow more cycles (and more relations) to
participate in $H_{k,w}Z$.  Whether this increases or decreases the size
of $H_{k,w}Z$ will just depend on the circumstances.   Sections 3 and 4
apply these concepts to $IC$ and toric varieties. This ends the summary.

The next topic is stratification independence.  It is important that the
$H_{k,w}Z$ not depend on the choice of the stratification.  This is
already known for mpih.  First suppose $w$ is zero, so the `where'
conditions are vacuous.  Suppose also that $\eta$ is a $k$-cycle, namely a
$\lambda$-coning of the global cycle $\eta_{(1)}$.  Already known is that
if the stratification is changed, an $\eta'$ equivalent to $\eta_{(1)}$
can be found.  This is a local result, and so the $\eta'$ can be found
close to $\eta_{(1)}$.  Thus, in this special case stratification
independence will follow if the $\lambda$-coning structure that converts
$\eta_{(1)}$ to $\eta$ can be extended to a small neighbourhood of $\eta$.
This is more a statement about the local conical structure of a stratified
topological space, than about the particular cycle $\eta$.

Now suppose $w$ is nonzero.  As before, $\eta'$ can be found close to
$\eta_{(1)}$.  The new difficulty is the $w(\eta_i)$.  By definition, this 
will cause no new problem, unless $w(\eta_i)$ is reduced.  In this case,
more must be done.  Previously, $\eta_{(1)}$ was moved
to $\eta'$, and it was assumed that $\eta'$ could be collapsed to the
$\eta_i$ used for $\eta$.  In this case, it is necessary to move the
$\eta_i$.  To do this, note that the collapsing $\eta \gg \eta_i$ will be
a permissable coning away, \emph{for some perversity other that the
middle}.  As stratification independence is known for all perversities,
this may provide a means of producing the $\eta'_i$.  As before, one hopes
that uniformity of the $\lambda$-coning will complete the proof.

These arguments indicate that stratification independence for local-global
homology will follow from uniformity of the $\lambda$-coning and some
modification of the existing methods.  There is however another approach. 
King \cite{bib.HK.TIIH} was able to prove stratification independence
directly, without use of sheaves in the derived category.  Habegger and
Saper \cite{bib.NH-LS.IC} found that this leads to an intersection
homology theory for generalised perversities and local systems.  Usually,
a local system on $Z$ is determined by its value on the generic stratum
$S_n\subset Z$. When generalised, information about behaviour at the
boundary (the smaller $S_i$) is recorded by the local system.  It may be
that local-global homology can be expressed as the homology of such a
local system. From this, stratification independence would follow.

(Certainly, the concepts are related.  If $\eta$ is a local-global cycle,
collapsing under $\lambda$ to $\eta_1$, and if $\eta_1\cap S_j$ is dense
in $\eta_1$, then the following holds.  There is a local system $L_j$ (in
the usual sense) on $S_j$ such that $\eta$ both determines and is
determined by a cycle $\eta'$ on $S_j$ with coefficients in $L_j$. 
(Because $S_j$ may not be closed, $\eta'$ might not be a compact cycle. 
This is a technical matter.)  Consider this $\eta'$.
The Goresky-MacPherson theory imposes
conditions on the dimension in which the closure $\bar\eta'$ of $\eta'$
meets the strata.  For local-global homology, $\eta'$ determines a
local-global cycle $\eta_0$ in a neighbourhood of $S_j$.  Whether or not
such an $\eta_0$ can be closed up to produce an $\eta$ is a delicate
matter, which depends not so much on how $\eta'$ meets the smaller strata,
but on how the local topology of $Z$ about these strata interacts with
$\eta_0$.  The generalised concept of a local system may provide a place
where this information can be stored, and used to control $\eta'$.)

This discussion has assumed that existing techniques, perhaps adapted,
will be applied to local-global homology.  However, stratification
independence is primarily a technical result on the local topology of
stratified spaces.  It may be that local-global homology will provide a
suitable language for describing this local topology.  One must show that
it causes no obstruction to the motion of cycles and boundaries, that is
required when the stratification is refined.  If this holds, then the
concept of local-global homology is already implicit in the proof of
stratification independence, and the consequences of King's paper become
less surprising.

Each Betti number of $\PDelta$ is the dimension of a vector space, and so
is nonnegative.  Thus, a linear function that computes such a Betti number
from $f\Delta$ is also a linear inequality on $f\Delta$, at least when
$\Delta$ has rational vertices.  Of special interest therefore are those
parts of homology theory, for which the Betti numbers are indeed linear
functions of $f\Delta$.  Ordinary homology does not have this
property~\cite{bib.McC.HTV}.  Middle perversity intersection homology does
\cite{bib.JD-FL.IHNP,bib.KF.IHTV}.  This is a consequence of deep results
in algebraic geometry, namely Deligne's proof~\cite{bib.NK.DPRH}
of the Weil conjectures and
the purity of mpih~\cite{bib.AB-JB-PD.FP}.  It is possible that at least
some part of the local-global homology theory will also have this purity
property.  An example of Bayer (personal communication) shows that certain
components of the extended $h$-vector are sometimes negative, and so are
not always Betti numbers.

(Bayer's example is $\Delta=BICCC\sqdot$, where $B$ is the bipyramid
operator, the combinatorial dual to the $I$ operator.  In $h\Delta$ the
term $xA\{1\}$ occurs with coefficient $-2$.  The interpretation of this
result requires some care.  It does not show that the whole of $h\Delta$
is unsuitable, or that local-global homology is a flawed concept.  The
formula for $h\Delta$ is the extrapolation to all polytopes of the
heuristically calculated formula for the various local-global homology
Betti numbers of the $IC$ varieties.  The component $xA\{1\}$ corresponds
to the gluing together along paths of local cycles, a construction that is
already known to yield a topological invariant.  More exactly, $xA\{1\}$
counts local $1$-cycles under this equivalence, on a $4$-dimensional variety
that has had its $0$-strata removed.  The corresponding Betti number is as
computed by $h\Delta$, for the $IC$ varieties, but clearly not in general.)

(Ordinary homology can be approached in the same way.  Use the simple case
formula $h(x,x+y)=f(x,y)$ to define a `pseudo $h$-vector' $\hquest\Delta$
for general polytopes.  The rules
\[
    \hquest (I\Delta) = (11) \hquest (\Delta) \> ;\quad
    \hquest (C\Delta) = (100\ldots0) + y\hquest (\Delta) \> ;
\]
give its transformation under $I$ and $C$.  The author suspects, but does
not known, that $\hquest\Delta$ gives the ordinary homology Betti numbers,
when $\Delta$ is an $IC$ variety.  However, when extrapolated to the
octahedron, $(1,-1,5,1)$ is the result.  Despite this negative value,
ordinary homology is still a topological invariant.  However, it is not
suitable for the study of general polytope flag vectors.)

For $\Delta$ rational the mpih part $(h_0,h_1,\ldots,h_n)$ of $h\Delta$
is not only nonnegative but also \emph{unimodal}.  This means that $h_0
\leq h_1 \leq \ldots$ and so on up to halfway.  It is a consequence of
strong Lefschetz.  The author suspects that the analogous result for
general polytopes will be as follows.  First compute $hC\Delta$.  This is
of course a linear function of $h\Delta$.  Now look at the coefficients of
$x^0y^0W$, where $W$ is a word in $A$ and $\{k\}$.  These numbers are
expected to be nonnegative.  This is an extension of the mpih result.  Of
course, algebraic geometry can prove such results only when $\Delta$ is
rational, for only then does $\PDelta$ exist.

The only method that will produce such results for general polytopes, that
the author can envision, is \emph{exact calculation of homology}. This
means constructing a complex of vector spaces (of length at most $n$),
whose Euler characteristic (alternating sum of dimensions) is the desired
Betti number.  Exact calculation then consists of proving that this
complex is exact, except possibly at one location.  The homology of this
complex at that location is then a vector space which as a consequence of
exactness has the desired dimension.  In \cite{bib.JF.CPLA} the author
introduces such complexes.  The proof of exactness is expected to be
difficult.  McMullen's proof \cite{bib.McM.SP} of strong Lefschetz for
simple polytopes, without the use of algebraic varieties, probably
contains a prototype of the arguments that will be needed.  There, the
Riemann-Hodge inequalities were vital in supporting the induction on the
dimension of the polytope.  Thus, one would like local-global homology to
support similar inequalities.

The process of constructing the complexes for exact calculation requires
that the linear function $h\Delta$ on polytope flag vectors be defined for
all flag vectors, whether of a polytope or not.  This allows one to talk
about the contribution made by an individual flag to $h\Delta$.  (One then
interprets this contribution as the dimension of a vector space associated
to the flag, and then assembles all these vector spaces into a complex. 
The boundary map is induced by the deletion of a single term from the
flag.  In \cite{bib.JF.CPLA}, the numerical contribution due to a flag is
derived from a study of the associated vector spaces, which is the
starting point.  This approach better respects the inner logic of the
exact calculation concept.)

The extension of $h\Delta$ to all flag vectors (not just those of
polytopes) can be done in the following way.  Suppose $\delta$ is an
$i$-face on $\Delta$.  The local combinatorial structure of $\Delta$ along
$\delta$ can be represented by an $(n-i-1)$ dimensional polytope, the
\emph{link $L_\delta$ of $\Delta$ along $\delta$}.  If $g_i$ is a linear
function on flag vectors then the expression
\begin{equation}
\label{eqn.h-sumg}
    h\Delta = \sum \nolimits_{ \delta \subseteq \Delta } g_i B
\end{equation}
is a linear function on the flag vector of $\Delta$.  (Here, $\delta$ runs
over all faces of $\Delta$, $i$ is the dimension of $\delta$, and $B$ is
the link along $\delta$.)  In this paragraph, $h\Delta$ is a linear
function that might or might not be the previously defined extended
$h$-vector.  Now use the rules
\begin{eqnarray}
    \label{eqn.g_0B}    g_0B &=& ChB -yhB \\
                        g_{i+1}B &=& yg_iB - g_i CB
\end{eqnarray}
and the initial value $g_0(\emptyset)=(1)$ to produce a recursive
definition of $g$ and $h$.  The $ChB$ in (\ref{eqn.g_0B}) stands for the
rule $\Ctilde$ of \S4, translated into a $x$, $y$, $A$ and $\{k\}$ rule,
and then applied to $hB$.  (Here, $B$ has dimension less than that of
$\Delta$, and so $hB$ is by induction already defined.)  This defines a
linear function $h\Delta$, which however can now be applied to
non-polytope flag vectors.  For polytope flag vectors, it agrees with the
previously defined value.  (Central to the proof of this is the expression
of the links on $I\Delta$ and $C\Delta$ in terms of those on $\Delta$. 
This approach does not rely on the $IC$ equation.)

The mpih part of (\ref{eqn.h-sumg}) is the usual formula 
\cite{bib.JD-FL.IHNP,bib.KF.IHTV,bib.RS.GHV} for these Betti numbers.  The
whole of (\ref{eqn.h-sumg}) can be `unwound' to express $h\Delta$ as a sum
of numerical contributions due to individual flags. The mpih part has been
presented in \cite{bib.MB.TASI}.  Note that just as the topological space
$\PDelta$ can be decomposed into cells in many ways, so its homology can
be computed in various ways.  To each suitable such method an extension of
$h\Delta$ to all flag vectors will follow, and vice versa.  An extension
which on simple polytopes reduces to the $h(x,x+y)=f(x,y)$ formula could
be very useful.  Finally, note that when $h\Delta$ is a sum of flag
contributions one can use `Morse theory' or shelling to compute
$h\Delta$.  Choose a linear `height' function, so that each vertex on
$\Delta$ has a distinct height.  Define the \emph{index} of a vertex $v$
to be the sum of the contributions due to the flag whose first term has
$v$ as its highest point.  It immediately follows that $h\Delta$ is the
sum, over all vertices, of their index.

Central to the problem of resolution of singularities (which is still open
in finite characteristic) is the discovery of suitable invariants of
singular varieties.  These should be well behaved under monoidal
transformation, and always permit the choice of a centre of
transformation, which will reduce the invariant.  When the invariant becomes
zero, the variety should be nonsingular.  It may be that some form of
local-global homology will have such properties.  One aspect of the
problem is this.  Suppose the singular locus consists of two lines meeting
at a point.  Along each line there is a locally constant singularity,
whose resolution is essentially a lower dimensional problem, which can by
induction be assumed solved.  The difficulty is at the meeting point of
the two lines.  Along each line there is a resolution process.  One needs
to know whether and when a transformation should be centered at the common
point.

In the study of this process one should use not only the local cycles due
to mpih, but also those due to the higher order local-global homology. 
(In addition, one should use integer rather than rational coefficients,
but more on this later.)  The local cycles generic along each of the two
lines will influence if not control the resolution process along that
line.  The concepts of local-global homology allow the interaction to be
studied at their common point, of these `controls' along the two lines of
the resolution process.  This argument indicates that local-global
homology contains the right sort of information, for it to act as a
suitable control on the resolution process.  It does not show that it
contains enough such information.

Torsion is when a cycle is not a boundary, but some multiple of it is. 
Using rational or real coefficients sets all such torsion cycles equal to
zero. This simplifies the theory, and for the study of Betti numbers and
the flag vectors of convex polytopes, such complications are not needed. 
For resolution of singularities, and more subtle geometric problems, the
situation is otherwise.  Here is an example.  Consider the affine surface
$X=\{xy=z^k\}$, for $k\geq 2$.  This has a singularity at the origin. 
However, there are no non-trivial local-global cycles on $X$, according to
the definitions of \S5.  The divisor class group ${\rm Cl}(X)$ however is
nontrivial, and all its elements are torsion.  For example, the line
$L=\{x=z=0\}$ is not defined by a principal ideal, whereas $kL$ is defined
by $\{x=0\}$.

Such information can be recorded by local-global homology, provided
integer coefficients and a different concept of a local cycle are used. It
is known that duality has a local expression, which pairs compact cycles
avoiding say the apex of a cone, and cycles of complementary dimension,
that meet the apex.  Duality, of course, ignores torsion cycles.  However,
non-trivial but torsion local cycles exist for $X=\{xy=z^k\}$, when the
complementary concept of cycle is used.  It is interesting that the
construction of local-global cycles associated to the formula
(\ref{eqn.h-sumg}) produces such cycles, rather than compact local-global
cycles.

It is natural to ask: does the vanishing of all such local-global cycles
(except mpih of course) imply that the variety is either nonsingular, or a
topological manifold?  Because the Brieskorn hypersurface singularity
$x_1^2 +x_2^2 +x_3^2 + y^3 + z^5=0$ is locally homeomorphic to the cone on
an exotic $7$-sphere \cite{bib.JM.MH7S}, which is knotted inside a $9$-sphere,
purely topological invariants are not enough to ensure the nonsingularity
of a variety $Z$.  The same may not be true for an embedded variety
$Z\subset \bfP_n$.

Ordinary homology and cohomology are functors.  This distinguishes them
from intersection homology.  For mpih and, one hopes, a good part of the
local-global homology, it is purity and formulae for Betti numbers that
are the characteristic properties.  Both concepts agree, of course, on
nonsingular varieties.  Thus, one of the geometric requirements on
local-global homology is as follows.  Suppose $f:Z_1 \to Z_2$ is a map
between, say, two projective varieties.  Suppose also that the strictly
local-global homology of both $Z_1$ and $Z_2$ vanishes.  It may be that
this condition ensures that the $Z_i$ are, say, rational homology
manifolds.  If this is so, then $f$ induces a map $f_*:H_\bullet Z_1 \to
H_\bullet Z_2$, which can be thought of as a point in $(H_\bullet
Z_1)^*\otimes H_\bullet Z_2$.  In the general case one would want $f$ to
induce some similar map or object in a space, which reduces to the
previous $f_*$ when the strictly local-global homology vanishes.

Similarly, when $Z$ is nonsingular its homology carries a ring structure. 
(It is this structure that supplies the pseudopower inequalities
\cite{bib.McM.NFSP} on the flag vector of a simple polytope.)  Elsewhere
\cite{bib.JF.IHRS}, the author provides for a general compact $Z$ a similar
structure, that in the nonsingular case reduces to the homology ring.

To close, some personal historical remarks, and acknowledgements.  In 1985, upon
reading the fundamental paper \cite{bib.MB-LB.gDS} of Bayer and Billera,
it became clear to the author that a full understanding of general convex
polytopes would require a far-reaching extension to the theory of
intersection homology.  The background to this insight came from Stanley's
proof \cite{bib.RS.NFSP} of the necessity of McMullen's conjectured
conditions \cite{bib.McM.NFSP} on the face vectors of simple polytopes, and
the proof by Billera and Lee of their sufficiency.  Danilov's exposition
\cite{bib.VD.GTV} of toric varieties, and McConnell's result \cite
{bib.McC.HTV}, also contributed.  It was also clear that once the extended
$h$-vector was known, the rest would soon follow.

Already in the simple case there is an interplay between the topological
definition of homology, $h$-vectors, and `combinatorial linear algebra' or
exact calculation.  The same holds for middle perversity intersection
homology.  In 1985 the Bernstein-Khovanskii-MacPherson formula was known,
although not published until 1987~\cite{bib.RS.GHV}.  Because of this
circle of ideas, knowledge of say the $h$-vector is sufficient in practice
to determine the other parts of the theory.  This led the author to find
the $IC$ equation \cite{bib.JF.MVIC} (again, known in 1985 but not
published until much later).  This moved the focus on to the rules for the
transformation of $h\Delta$ under $I$ and $C$.

In the early years of the search for these rule, two related and
erroneous ideas were influential.  The first is that the rule for $I$ (as
via K\"unneth) should be multiplication by $(x+y)$.  The second is that
the generalised Dehn-Somerville equations should be expressed by
$h\Delta$ being palindromic.  (Also unhelpful was an undue concentration on
the formula for $h\Delta$.)  The crucial step that lead to these
assumptions being dropped took place in 1993.  Loosely speaking, it was
the discovery of special cases of the `gluing local cycles together along
paths' construction.  At that time how to form families of local cycles,
or in other words the `local-global' concept, was still mysterious.

In late 1995 the present formula for $h\Delta$ was discovered, as a
solution to the various geometric, topological and combinatorial
constraints that were known.  It was not at that time properly understood
by the author.  This definition was in 1996 pushed around the circle of
ideas, to produce first a combinatorial linear algebra construction for
$H_\bullet\Delta$, and then the topological definition.  That all this can
be done indicates that the definition of $h\Delta$ is correct.  Once these
advances had been understood, it was then possible to return to the
derivation of the formula for $h\Delta$, and put it on a proper footing. 
It was only at this point that the fundamental concepts became clear.  In
1997, the analogue to the ring structure was found.  Its relation to the
concepts presented here is, at the time of writing, still under
investigation.

Many of the results in this paper were first made available as preprints
and the like, in 1996 and 1997.  The difficulties encountered by the
readers have lead to many revisions in the exposition, and clarification
of the basic concepts, both on paper and in the author's understanding.
Part of the difficulty is the extent of the circle of ideas, which passes
through several areas of mathematics.  Any one of topology, combinatorics,
linear algebra and intersection theory can be chosen as the starting
point.  Another difficulty is that much of the intuition and guidance
comes from perhaps uncomplicated examples and points of view that have, by
and large, not yet been put into print.

Finally, thanks are due to Marge Bayer, Lou Billera, Martin Hyland, Gil
Kalai, Frances Kirwan, Carl Lee, Robert MacPherson, Mark McConnell, Peter
McMullen, Rick Scott, Richard Stanley and David Yavin variously for their
tolerant interest in the author's previous attempts to find $h\Delta$, to
define $H_{k,w}Z$, and to express these results.

\end{document}